\newcommand{\les}{\lesssim}

\def\piT{{\,^{(\T)}\pi}}
\def\Qr{\mbox{Qr}}

\def\Null{\dot{\NN}^{-}}

\def\GG{{\cal G}}
\def\OO{{\cal O}}
\def\QQ{{\cal Q}}

\def\VV{{\cal V}}
\def\TT{{\cal T}}

\def\P{{\bf P}}
\def\vphi{\varphi}
\documentstyle[amssymb,amsfonts]{amsart}

\newenvironment{proof}{\noindent {\bf Proof} }{\endprf\par}
\def \endprf{\hfill  {\vrule height6pt width6pt depth0pt}\medskip}
\def\emph#1{{\it #1}}
\def\textbf#1{{\bf #1}}
\newcommand{\bea}{\begin{eqnarray}} 
\newcommand{\eea}{\end{eqnarray}}
\def\beaa{\begin{eqnarray*}}
\def\eeaa{\end{eqnarray*}}

\def\ba{\begin{array}}
\def\ea{\end{array}}
\def\be#1{\begin{equation} \label{#1}}
\def \eeq{\end{equation}}

\newcommand{\nn}{\nonumber}

\def\a{{\alpha}}
\def\b{{\beta}}
\def\ga{\gamma}

\def\de{\delta}
\def\De{\Delta}
\def\ep{\epsilon}
\def\eps{\epsilon}
\def\la{\lambda}
\def\La{\Lambda}
\def\si{\sigma}
\def\Si{\Sigma}
\def\om{\omega}
\def\Om{\Omega}

\def\ze{\zeta}
\def\nab{\nabla}
\def\varep{\varepsilon}

\def\aa{{\underline{\a}}}
\def\bb{{\underline{\b}}}
\def\bb{{\underline{\b}}}

\def\Lb{{\underline{L}}}

\newcommand{\chih}{\hat{\chi}}
\newcommand{\chib}{\underline{\chi}}

\def\chih{\hat{\chi}}
\def\trch{\mbox{tr}\chi}
\def\tr{\mbox{tr}}

\def\Xb{{\underline X}}
\def\Yb{{\underline Y}}

\def\MM{{\cal M}}
\def\NN{{\cal N}}
\def\II{{\cal I}}
\def\FF{{\cal F}}

\def\JJ{{\cal J}}
\def\KK{{\cal K}}
\def\Lie{{\cal L}}
\def\DD{{\cal D}}

\def\RR{{\cal R}}

\def\D{{\bf D}}
\def\F{{\bf F}}

\def\I{{\bf I}}

\def\M{{\bf M}}
\def\L{{\bf L}}

\def\Q{{\bf Q}}
\def\R{{\bf R}}

\def\T{{\bf T}}
\def\g{{\bf g}}
\def\m{{\bf m}}

\def\LLb{{\bf \Lb}}

\def\RRR{{\Bbb R}}

\def\SSS{{\Bbb S}}

\def\nabb{\overline{\nab}}

\def\expp{\mbox{exp}}

\def\pr{\partial}
\newcommand{\ddd}{\nab}

\def\c{\cdot}

\def\f14{{\frac{1}{4}}}
\def\f12{{\frac{1}{2}}}

\def\En{\mbox{En}}

\def\dual{{\,^\star \mkern-3mu}}
\def\2{{\overline 2}}

\begin{document}
\theoremstyle{plain}
  \newtheorem{theorem}[subsection]{Theorem}
  \newtheorem{conjecture}[subsection]{Conjecture}
  \newtheorem{proposition}[subsection]{Proposition}
  \newtheorem{lemma}[subsection]{Lemma}
  \newtheorem{corollary}[subsection]{Corollary}

\theoremstyle{remark}
  \newtheorem{remark}[subsection]{Remark}
  \newtheorem{remarks}[subsection]{Remarks}

\theoremstyle{definition}
  \newtheorem{definition}[subsection]{Definition}

\include{psfig}
\title[Curvature ]{On the radius of injectivity of null hypersurfaces.}
\author{Sergiu Klainerman}
\address{Department of Mathematics, Princeton University,
 Princeton NJ 08544}
\email{ seri@@math.princeton.edu}

\author{Igor Rodnianski}
\address{Department of Mathematics, Princeton University, 
Princeton NJ 08544}
\email{ irod@@math.princeton.edu}
\subjclass{35J10\newline\newline
The first author is partially supported by NSF grant 
DMS-0070696. The second author is partially 
supported by NSF grant DMS-0406627. Part of this work was done while he was 
visiting Department of Mathematics at MIT
}
\maketitle
\vspace{-0.3in}
\section{Introduction} This paper is concerned  with the
regularity properties of boundaries $\NN^{-}(p)=\pr\II^{-}(p)$ 
 of  pasts (future) of points in a  $3+1$ Lorentzian manifold $(\M, \g)$. 
The  past of a point  $p$, denoted $\II^{-}(p)$, is the collection
of points that can be reached by a past directed  
time-like curve from $p$. As it is well known 
the past boundaries $\NN^{-}(p)$ play a crucial role 
in understanding the causal structure of Lorentzian manifolds
and the propagation properties of linear and nonlinear waves, e.g.
in flat space-time the null cone $\NN^{-}(p)$ is exactly the propagation
set of solutions to the standard wave equation  with a 
Dirac measure  source point at $p$.  However these past boundaries
 fail, in general,  to be smooth even in a  smooth, curved, lorentzian 
space-time; one can only guarantee that $\NN^{-}(p)$ is
a Lipschitz, achronal,  $3$-dimensional manifold without boundary ruled
by  inextendible null geodesics from $p$,
see \cite{HE}. In fact $\NN^{-}(p)\setminus\{p\}$ is smooth in a small
neighborhood of $p$ but fails to be so in the large 
 because of conjugate points, resulting in  formation of caustics,
or because of  intersections of distinct null geodesics from $p$.
Providing a lower bound for the radius of injectivity of the
sets  $\NN^{-}(p)$ 
 is  thus  an essential step in understanding the more
refined  properties
of solutions to linear and nonlinear wave equations a Lorentzian background.

 The  phenomenon described above  is also present in Riemannian geometry in
connection to geodesic coordinates relative to a point, yet in that case the
presence of conjugate  or cut-locus points has nothing to do
with the regularity of the manifold itself.  In that sense   
lower bounds  for  the radius of injectivity of a  Riemannian manifold  are 
important  only   in  so far as   geodesic normal coordinates,  and their applications\footnote{such as, for example,  Sobolev inequalities   or  the 
finiteness theorem of Cheeger.},   are concerned.
 Thus, for example, lower bounds for the  radius of injectivity
can  sometimes be replaced by lower bounds for the   harmonic radius, which plays
an important role in  Cheeger-Gromov theory, see e.g. \cite{A2}.

In this paper we investigate regularity of past boundaries $\NN^-(p)$ 
in Einstein vacuum space-times, i.e., Lorentzian manifolds $(\M,\g)$
with the Ricci flat metric,
$
\R_{\a\b}(\g)=0$.
We provide
 conditions on  an Einstein vacuum space-time $(\M, \g)$,
compatible with bounded $L^2$ curvature,
 which are sufficient to ensure the 
local  non-degeneracy of  $\NN^-(p)$.  
More precisely  we provide  a  uniform lower bound
on the radius of injectivity of  the null boundaries $\NN^{-}(p)$
of the causal past sets  $\JJ^{-}(p)$ in terms of the Riemann curvature
flux on $\NN^-(p)$ and some other natural assumptions\footnote{arising
specifically in applications to the the problem of a break-down criteria
in General Relativity discussed in \cite{KR5}.} on
$({\M},\g)$. Such lower bounds are essential
in understanding the causal structure and  the related
 propagation properties of  solutions to the Einstein equations. They are
particularly important in construction of an effective  Kirchoff-Sobolev type
parametrix for solutions of wave equations on $\M$,
see \cite{KR4}. Such parametrices are used in \cite{KR5}
to prove a large data break-down criterion for solutions
of the Einstein-vacuum equations.

This  work complements our series 
of papers \cite{KR1}-\cite{KR3}. The methods
of  \cite{KR1}-\cite{KR3} can be  adapted\footnote{In 
\cite{KR1}-\cite{KR3}  we have considered the case 
of the congruence of outgoing future null geodesics initiating
on a $2$-surface  $S_0$ embedded in a space-like hypersurface $\Si_0$.
The extension of our results to null cones from a point forms
the subject of Qian Wang's Princeton 2006 PhD thesis.}
to prove 
  lower bounds on  the geodesic radius of  conjugacy of the 
congruence of past null geodesics from $p$ which depends
only on the  geodesic (reduced)  flux of curvature, i.e.  
  an $L^2$ integral norm  along $\NN^{-}(p)$  of  tangential components
of the  Riemann curvature tensor $\R=\R(\g)$,
see  section \ref{sect:geod-curv-flux} for a precise definition.  
 It is however possible that  the radius of conjugacy
of the null congruence is bounded from below and yet there are
past null geodesics
 form a point $p$ intersecting again at points arbitrarily close
to $p$. Indeed, this can happen on a flat Lorentzian  manifold 
such as $\M={\Bbb T}^3\times \RRR$ where ${\Bbb T}^3$ is the torus obtained
by identifying the opposite sides of a lattice of period $L$
and metric induced  by the standard Minkowski metric. Clearly
there can be no conjugate points for the  congruence of past
or future null geodesics from a point and yet there are plenty of
 distinct null geodesics from a point $p$ in $\M$  which  intersect 
on a time scale proportional to $L$. There can be thus no lower bounds
on the null radius of injectivity expressed only in terms of bounds for the
curvature tensor $\R$. This problem occurs, of course, also in Riemannian
geometry where we can control the  radius of conjugacy in terms 
of uniform bounds  for the curvature tensor, yet, in order to control
the   radius of injectivity we need to make other 
  geometric assumptions such as, in the case of compact Riemannian
manifolds,   lower bounds on   volume  and upper bounds for
its diameter.  It should thus come
as no surprise that we also need, in addition to 
bounds for the curvature flux,  other  assumptions on
 the geometry of    solutions to the Einstein equations
in order to ensure  control on  the null cut-locus of points
in $\M$ 
 and  obtain
 lower bounds for  the  null
radius of injectivity. 

In this paper  we give  sufficient conditions,
expressed relative  to a space-like foliation $\Si_t$ given by the level
surfaces of a regular time function $t$ with unit future
 normal $\T$.  We discuss two related  results. Both
are based on a space-time assumption on the uniform boundedness
of the deformation tensor $\piT=\Lie_\T\g$  and boundedness of the
 $L^2$  norm of the   curvature tensor on a fixed  slice $\Si_0=\Si_{t_0}$
of the foliation. Standard energy estimates,
based on the Bel-Robinson tensor, allows to get a uniform
control on the $L^2$ norm of the curvature tensor 
on all slices.   In the first result we also assume 
that  every point  of the space-time  admits a sufficiently large coordinate 
patch with a  system of coordinates in which the
 Lorentz metric $\g$ is  close   to  a flat
 Minkowski  metric. In the second result we  dispense of 
the latter  condition by showing how such coordinates  
can be constructed, dynamically, from  a given
coordinate system on the initial slice $\Si_0$.
 Though the second result  (Main Theorem II), is more appropriate 
for applications, the  main new ideas of the paper 
appear in section 2  related to the proof
of the first result (Main Theorem I).

 The  energy estimates  mentioned above
also provide   uniform control 
on the  geodesic (reduced) curvature flux along 
 the null boundaries $\NN^{-}(p)$ and thus,
according to \cite{KR1}-\cite{KR3}, give
control on the radius of conjugacy of the corresponding
null congruence. There is however an important  subtlety involved here.
Past points of intersection, distinct   null geodesics   
from $p$ are no longer on the boundary of $\JJ^{-}(p)$ 
and therefore  the energy estimates mentioned above do
not apply.  Consequently  we cannot simply 
apply the results \cite{KR1}-\cite{KR3} and estimate
the  null radius of conjugacy independent of the cut-locus,
but  have to treat them together by a delicate  boot-strap
argument. The main new ideas of this paper concern estimates
for the cut locus, i.e. establishing  lower bounds, with
respect to the time parameter $t$,  for the
points of intersection of   distinct  past null geodesics
 from $p$.
Though our results,  as  formulated, hold only for  
Einstein vacuum manifolds  our method of proof  in section 2
 can  be extended to general Lorentzian manifolds if 
we make, in addition to the assumptions
 mentioned above,  uniform norm assumptions for the 
curvature tensor $\R$. Our results seem to be new even
in this vastly simplified  case, indeed we are not
aware of any non-trivial results concerning the null
radius of injectivity for  Lorentzian manifolds.

We now want to make a comparison with the corresponding 
picture in Riemannian geometry. In general, all known lower 
bounds on the radius of injectivity require some pointwise control
of the curvature. 
The gold standard in this regard is a theorem of Cheeger providing a lower bound 
on the radius of injectivity in terms of pointwise bounds on the sectional
curvature and diameter, and a lower bound on the volume of a compact manifold, 
see \cite{Ch}. 
Similar to our case, the problem in Riemannian geometry  splits into a lower
bound on the radius of conjugacy and an estimate on the length of
the shortest geodesic loop. The radius of conjugacy is intimately 
tied to  point-wise bounds on curvature, via the Jacobi equation. 
The estimate for the length of the shortest geodesic loop  relied, traditionally,
 on the Toponogov's Theorem, which again needs pointwise bounds
on the curvature. These two problems can
be naturally separated in Riemannian geometry\footnote{something that we do not know
how to do in the Lorentzian case} and while the radius of conjugacy 
requires pointwise assumptions on curvature, a
lower bound on the length of the shortest geodesic can be given  under  weaker,  integral, 
assumptions on curvature. The best result in the latter direction, to our knowledge, 
is due to Petersen-Steingold-Wei, which  in addition to the usual diameter and volume 
conditions on an $n$-dimensional compact manifold, 
requires {\it smallness} of the $L^p$-norm of sectional curvature, for
$p>n-1$, \cite{PSW}.
 Once more, we want to re-emphasize the fact that in Riemannian geometry
 lower  bounds on the radius of injectivity require  pointwise bounds
for the curvature, yet this restriction can be often overcome in applications 
by replacing it with bounds for  more  flexible  geometric quantities.
   A case in point is the  Anderson-Cheeger
result \cite{AC} which proves a finiteness theorem    under
   pointwise assumptions  on the Ricci curvature and $L^{\frac n2}$ bounds on the
 full Riemann curvature  tensor. Unlike the classical result of Cheeger,
see \cite{Ch}, the radius of injectivity need not, and cannot in general,
 be  estimated.

We should note that the works in Riemannian geometry,
cited above, have been largely stimulated
by applications to the Cheeger-Gromov theory. Applications of this theory to 
General Relativity have been pioneered by M. Anderson, \cite{A2} and \cite{A3}, see also 
\cite{KR5} for further applications. M. Anderson has, in particular, been
interested  in the  possibility of  transferring some aspects of Cheeger-Gromov theory 
to the Lorentzian setting. With this in mind he   
 has  proved existence of a  special 
space-time coordinate system for  Einstein  vacuum space-times, under pointwise
assumptions on the space-time Riemann curvature tensor, see  \cite{A1}.
 Another example  of a  global  Riemannian geometric  which  has been
successfully transplanted to Lorentzian setting is Galloway's null splitting theorem,
see  \cite{G}. 

We note that,   in the Riemannian setting,  the radius of injectivity 
and shortest geodesic loop estimates depend  crucially  on  lower bounds for the
volume of  the   manifold, as confirmed by the  example of a thin flat torus. By contrast, 
 the notion of volume  of a   null hypersurface in Lorentzian geometry  
is not well-defined,  as the restriction of the space-time metric to a null hypersurface 
is degenerate. We are forced  to replace the condition on the volume of $\NN^-(p)$ 
with the condition on the volume of the $3$-dimensional domain obtained by intersecting 
the causal past of $p$ with the level hypersurfaces of a time function $t$. To be more
precise our assumption on existence of a coordinate system in which the metric 
$\g$ is close to the Minkowski metric will allow us to prove that the volume of these
domains, at time $t<t(p)$, are  close to the volume of the Euclidean
 ball of radius $t(p)-t$,
where
$t(p)$ denotes the  value of the parameter  $t$ at $p$.
 Another ingredient of our proof is an argument showing that, at the 
first time of intersection $q$,  past null geodesics from a point $p$ 
meet each other at angle $\pi$,
viewed with respect to  the tangent space of the  space-like hypersurface 
$t=t(q)$. We also
show that we can find a point $p$, such that  the above property  holds 
both at $p$ and the first intersection point $q$.
 Finally, we give
 a geometric comparison
argument showing that an existence of a pair of null geodesics from a point $p$ 
meeting each other at angle $\pi$ both at $p$ and at the point of first intersection 
violates the structure of the past set $\JJ^-(p)$, if the time of
 intersection is too close to the value $t(p)$.

\section{Main results}

 We consider a space-time $(\M, \g)$ verifying the Einstein -vacuum
equations,
\be{eq:Einst-vacuum}
\R_{\a\b} =0.
\end{equation}
and assume that a part of space-time $\MM_I\subset \M$ is foliated 
by the level hypersurfaces of a time function $t$, monotonically
increasing towards future in the interval $I\subset \RRR$. Without loss
of generality we shall assume that the length of $I$, verifies,
\beaa
|I|\ge 1.
\eeaa
Let $\Si_0$ be a fixed
leaf  of the $t$ foliation.  Starting with a local coordinates chart $U\subset \Si_0$
and coordinates  $(x^1,x^2, x^3)$  we parametrize 
the domain  $I\times U\subset {\cal M}_I$ 
 with  transported  coordinates  $(t,x^1,x^2, x^3)$ 
obtained by following the integral curves of
$\T$, the future unit normal to $\Si_t$.
 The space-time metric $\g$ on $I\times U$ then takes the  form 
\be{eq:g-transported}
\g=-n^2 dt^2 + g_{ij} dx^i dx^j,
\end{equation}
where $n$ is the lapse function of the $t$ foliation and $g$ is the restriction of 
the metric $\g$ to the surfaces $\Si_t$ of constant $t$.  We shall assume
that the space-time region $\MM_I$ is globally hyperbolic, i.e. every
 causal curve from a point $p\in \MM_I$  intersects  $\Si_0$ at
 precisely one point.

  The
second fundamental  form of $\Si_t$ is defined by,
$$
k(X,Y)=\g(\D_X \T,Y),\quad \forall X,Y\in T(\Si_t),
$$
with $\D$  the Levi-Civita  covariant derivative.  Observe that,
\be{eq:t-sec-fund}
\pr_t g_{ij} = -\frac 12 n\, k_{ij}.
\end{equation}
   We denote by $\tr k$ the trace  of
$k$ relative to $g$, i.e. $\tr k = g^{ij} k_{ij}$.     We also assume 
the surfaces $\Si_t$ are either  compact or asymptotically flat.

Given a unit time-like normal $\T$ we can define a pointwise norm $|\Pi(p)|$ of any 
space-time tensor $\Pi$ with the help of the decomposition 
$$
X= X^0 \T + \underline X,\quad X\in T{\cal M},\quad \underline X\in T\Si_t
$$
The norm $|\Pi(p)|$ is then defined relative to the Riemannian metric,
\be{eq:riem-metric}
\bar{\g}(X,Y)=X^0\c Y^0+g(\Xb, \Yb).
\end{equation}
 We denote by $\|\Pi(t)\|_{L^p}$ the $L^p$ norm of $\Pi$
on $\Si_t$. More precisely,
\beaa
\|\Pi(t)\|_{L^p}=\int_{\Si_t}|\Pi|^p dv_g
\eeaa
with $dv_g$ the volume element of the metric $g$ of $\Si_t$.

Let $\piT=\Lie_\T\g$ be the deformation tensor of $\T$. The components of 
$ \piT$ are given by 
$$
\piT_{00}=0,\quad \piT_{0i}=\nab_i\log n,\quad \piT_{ij}=-2k_{ij}
$$
\subsection{Main assumptions}
We make the space-time assumption,
\bea
N_0^{-1}&\le& n\,\,\le\, N_0\label{eq:assumption-n}\\
|I|\c \sup_{t\in I}\|\pi(t)\|_{L^\infty} &\le& \KK_0<\infty.\label{eq:assumption-T}
\eea
where $N_0, \KK_0>0$ are given numbers and  $|I|$ denotes the length of the time
 interval $I\subset\RRR$.
We also make the following assumptions on the initial  hypersurface $\Si_0$,

{\bf I1.} \quad There exists a  covering 
of $\Si_0$  by    charts $U$ such
that for any fixed chart, the induced metric $g$  verifies
\be{eq:assumption-I1}
 I_0^{-1}|\xi|^2\le g_{ij}(x) \xi_i\xi_j \le I_0|\xi|^2, \qquad \forall
x\in U
\end{equation}
with $I_0$ a fixed positive number.
Moreover there
 exists a number $\rho_0>0$ such that every point $y\in \Si_0$  admits a neighborhood
$B$, included in a neighborhood chart   $U$, such that $B$ is precisely the Euclidean 
ball $B=B_{\rho_0}^{(e)}(y)$ relative to the local  coordinates  in $U$.

{\bf I2.}\quad The  Ricci  curvature of the initial slice $\Si_0$
verifies,
\be{eq:initial-L2-curv}
\|\R\|_{L^2(\Si_0)}\le \RR_0<\infty.
\end{equation}
{\bf Remark.}\quad If $\Si_0$ is compact the existence of $\rho_0>0$ 
is guaranteed by the existence of the coordinates charts $U$ verifying 
\eqref{eq:assumption-I1}. More precisely we have:
\begin{lemma} \label{le:initial-rho} If $\Si_0$ is compact and 
has a system of coordinate charts $U$ verifying \eqref{eq:assumption-I1},
 there
must exist a number $\rho_0>0$ such that every point $y\in \Si_0$  admits a neighborhood
$B$, included in a neighborhood chart   $U$, such that $B$ is precisely the Euclidean 
ball $B=B_{\rho_0}^{(e)}(y)$ relative to the local  coordinates  in $U$.
\end{lemma}
\begin{proof}:\quad 
According to our assumption  every point $x\in \Si_0$ 
belongs to a coordinate patch $U$. Let $r(x)>0$ be such
that the euclidean ball, with respect to the coordinates 
of $U$,  centered at $x$ of radius $r(x)$ is included in $U$.
Due to compactness of $\Si_0$ we can find a finite number 
of points $x_1,\ldots x_N$ such that the balls $B^{(e)}_{r_j/2}(x_j)$,  with $r_j=r(x_j)$
 for $j=1,\ldots ,N$,
cover $\Si_0$. Thus any $y\in \Si_0$ must belong to
 a ball  $B^{(e)}_{r_j/2}(x_j)\subset B^{(e)}_{r_j}(x_j)\subset U$, for some $U$. 
Therefore  the ball $B^{(e)}_{r_j/2}(y)\subset U$.   We then choose
 $\rho_0=\min_{j=1}^Nr_j/2$.
\end{proof}

 \subsection{ Null boundaries of $\JJ^{-}(p)$ }Starting with any point $p \in \MM_I\subset\M$,
 we denote by $\JJ^{-}(p)$ the causal past of $p$,   by $\II^{-}(p)$ 
its  interior and by $\NN^{-}(p)$ its null boundary all restricted to 
the region $\MM_I$ under consideration.

 In general $\NN^{-}(p)$ is
an achronal, Lipschitz  hypersurface, ruled
by the set of past  null geodesics  from $p$. 
We parametrize these geodesics   with respect to the  future, unit, time-like vector
 $\T_p$. Then,  for every direction $\omega\in \SSS^2$,
 with  $\SSS^2$ denoting the standard sphere in $\RRR^3$,  consider the  
 null vector $\ell_\omega$ in $T_p\M$,  
\be{eq:norm-geodd}
\g(\ell_\omega, \T_p)=1,
\end{equation}
  and associate to it the  past null
geodesic $\ga_\omega(s)$ with initial data $\ga_\omega(0)=p$ and 
$\dot \ga_\omega(0)=\ell_\omega$. 
We further define a null vectorfield $\L$ on $\NN^-(p)$ according to 
$$
\L(\gamma_\omega(s))=\dot\ga_\omega(s).
$$
 $\L$ may only be  smooth almost everywhere on $\NN^{-}(p)$ and 
can be multi-valued on a set of exceptional points. 
We can choose the parameter $s$ in 
such a way so that $\L=\dot\ga_\om(s)$ is geodesic and $\L(s)=1$.

 For a sufficiently small $\de>0$ the  exponential map  $\GG$  defined by, 
\be{eq:exp-map}
\GG:\,\,(s,\om)\to \ga_\om(s)
\end{equation}
is a  diffeomorphism from $[0,\de)\times \SSS^2$
to its image in $\NN^{-}(p)$.
 Moreover  for each $\om\in \SSS^2$ either $\ga_\om (s)$ can be
 continued for all positive  values of $s$ or there exists a   value
$ s_*(\om)$ beyond which  the points $\ga_\om(s)$ are no longer on the boundary
 $\NN^{-}(p)$ of $\JJ^{-}(p)$ but rather in its interior, see \cite{HE}.  We call
such points terminal points of $\NN^{-}(p)$. We say that a terminal point $q=\ga_\om(s_*)$ 
is a conjugate terminal point  if the map  $\GG$ is singular at $(s_{*},\om)$. 
A terminal point  $q=\ga_\om(s_*)$ is said to be a cut locus terminal point if  
the map $\GG$ is nonsingular at  $(s_{*},\om)$ and there exists another 
null geodesic from $p$, passing through $q$.  

 Thus    $\NN^{-}(p)$ 
 is a smooth manifold at all points except the vertex $p$ and
 the  terminal points of its past null geodesic generators. We denote
by $\TT^{-}(p)$ the set of all terminal points and by 
$\Null(p)=\NN^{-}(p)\setminus \TT^-(p)$ the  smooth  portion of $\NN^-(p)$.
The set $\GG^{-1}(\TT^{-}(p))$
has measure zero relative to the standard  measure $ds da_{\SSS^2}$ of the cone
$[0, \infty)\times \SSS^2$. We will denote by $dA_{\NN^{-}(p)}$ the 
corresponding    measure on $\NN^{-}(p)$. Observe that the definition is not intrinsic,
it depends in fact   on the  normalization condition \eqref{eq:norm-geodd}.

\begin{definition}\label{def:rad-inj}We define $i_*(p)$  to be the supremum over
 all the values $s>0$ for which the exponential map 
$\GG:(s,\om)\to \ga_\om(s)$ is a global diffeomorphism.
We shall refer to $i_*(p)$ as the null  radius of injectivity 
at  $p$ relative to the geodesic foliation 
defined by \eqref{eq:norm-geodd}.
\end{definition}
 {\bf Remark.}\quad  Unlike in Riemannian geometry
where the radius of injectivity
is defined with respect to the distance function,
the definition above depends on the  normalization
 \eqref{eq:norm-geodd}.

 \begin{definition}
We denote by $\ell_*(p)$ the smallest value of $s$ for which
there exist two distinct  null geodesics
 $\ga_{\om_1}(s), \ga_{\om_2}(s)$
 from $p$ which intersect  at a point for which 
the  corresponding value (smallest  for $\ga_{\om_1}$ and $\ga_{\om_2}$) of the affine parameter 
 is    $s=\ell_*(p)$.
\end{definition}

\begin{definition}
 Let $s_*(p)$ denote the supremum over all values of $s>0$ 
such that  the exponential map is
a local diffeomorphism  on $[0,s)\times\SSS^2$ . We shall refer to 
$s_*(p)$ as the  null radius of conjugacy of the point $p$.  
\end{definition}
Clearly,
\be{eq:ils}
i_*(p)= \min(l_{*}(p), s_*(p)) 
\end{equation}
  The first  goal of this paper is to prove  the following theorem
 concerning a lower bound for the radius of
 injectivity $i_{*}$of a space-time region $\MM_I$, under the following 
assumption:

{\bf Assumption  C. }\quad {\it  Every  point $p\in \MM_I$ admits a 
coordinate neighborhood 
$I_p\times U_p$ such that 
$U_p$ contains a geodesic ball $B_{r_0}(p)$ and 
\be{eq:const-r0}
\sup_{t,t'\in I_p}|t-t'|\ge r_0. 
\end{equation}
We assume that on $I_p\times U_p$ there exists a system of transported coordinates
\eqref{eq:g-transported}  close to the flat Minkowski metric $-n(p)^2dt^2+\de_{ij} dx^i
dx^j$. More precisely,
\be{eq:n-1}
|n(t,x)-n(p)|<\ep
\end{equation}
\be{eq:g-de}
 |g_{ij}(t,x)-\de_{ij}|<\ep
\end{equation}
 where $n(p)$ denotes the value of the lapse  $n$ at $p$.
 }
\begin{theorem}[Main Theorem I] \label{thm:main1} Assume that $\MM_I$ is globally hyperbolic
and verifies 
 the main assumptions \eqref{eq:assumption-n}, \eqref{eq:assumption-T},
\eqref{eq:initial-L2-curv} as well as assumption {\bf C} above.
We also assume that $\MM_I$
contains a future, compact set
${\cal D}\subset\MM_I$  such that there exists a positive constant $\de_0$
for any point $q\in {\cal D}^c$ we have $\ell_*(q)>\de_0$.

Then, for sufficiently  small $\ep>0$, there exists 
 a positive number
$i_{*}>0$, depending only on $\de_0, r_0$, $N_0$,  $\KK_0$ and  $\RR_0$, 
such that, for all $p\in \MM_I$, 
\be{eq:lower-bound-inj}
i_{*}(p)>i_*
\end{equation}
\end{theorem}
 Assumption {\bf C} of  theorem \ref{thm:main1} 
can in fact be eliminated according to the following. 
\begin{theorem}[ Main Theorem II]\label{thm:main2} Assume that $\MM_I$ is globally hyperbolic
and  verifies the assumptions \eqref{eq:assumption-n}, \eqref{eq:assumption-T} as well
as the initial assumptions {\bf I1} and {\bf I2}.  Assume  also that  $\MM_I$
contains a future, compact set
${\cal D}\subset\MM_I$  such that $\ell_*(q)>\de_0$ for any point $q\in {\cal D}^c$,
 for some  $\de_0>0$.

 There exists 
 a positive number
$i_{*}>0$, depending only on  $I_0, \mu_0$, $\de_0$,  $N_0$,  $\KK_0$ and  $\RR_0$, 
such that, for all $p\in \MM_I$, 
\be{eq:lower-bound-inj}
i_{*}(p)>i_*
\end{equation}
\end{theorem}
{\bf Remark.}\quad 
Observe  that the last assumption of both theorems, concerning
a lower bound for $l_*$ outside  a sufficiently large future compact set,
is superfluous on a manifold with compact initial slice\footnote{A similar statement
can be made if $\Si_0$ is  asymptotically  flat.} $\Si_0$. Thus, for manifolds
$\MM_I=I\times\Si_0$, with $\Si_0$ compact $i_*$ depends only on the constants $I_0$,$N_0$,
$\KK_0$ and $\RR_0$.

The first key step in the proof of the Main Theorem is a lower bound on
the radius of conjugacy $s_*(p)$. 
\begin{theorem}\label{thm:main-conditional}
 There exists
a sufficiently small  constant $\de_*>0$, depending only on 
$ \KK_0$ and $\RR_0$ such that,
 for any $p\in \MM_I$
we must have,
$$s_*(p) >\min (\ell_*(p),\de_*).$$
\end{theorem}
The proof of Theorem \ref{thm:main-conditional} crucially relies on the 
results obtained in \cite{KR1}-\cite{KR3}\footnote{and an extension in Q. Wang's 
thesis, Princeton University, 2006.}. The discussion of these results and 
their reduction to Theorem \ref{thm:main-conditional} will be discussed in Section 
\ref{Conjugacy}.

\subsection{Connection to  $t$-foliation}We reinterpret the result
of theorem  \ref{thm:main-conditional} relative to the $\Si_t$ foliation.
For this we first make the following definition.
\begin{definition}\label{def:rad-inj-t} Given $p\in \MM_I$ we define $i_*(p, t)$ 
 to be the
 supremum  over all the values $t(p)-t$ for which the exponential map 
$$\GG:(t,\om)\to \ga_\om(t)=\ga_\om(s(t))$$ is a global diffeomorphism.
We shall refer to $i_*(p, t)$ as the null  radius of injectivity 
at  $p$ relative to the $t$-foliation. We denote by $\ell_*(p,t)$ the
supremum  over all the values $t(p)-t, \, t<t(p),$  for
which there exist two distinct  null geodesics
 $\ga_{\om_1}, \ga_{\om_2}$,
 from $p$ which intersect  at a point on $\Si_t$. Similarly, we let
 $s_*(p,t)$ be the supremum of $t(p)-t$ for which the exponential 
 map $\GG(t,\omega)$ is a local diffeomorphism.
\end{definition}
The results leading up to the proof of Theorem \ref{thm:main-conditional} 
also imply the following
\begin{theorem}\label{thm:main-conditional-t}
 There exists
a sufficiently small  constant $\de_*>0$, depending only on $N_0$,
$ \KK_0$ and $\RR_0$ such that,
 for any $p\in \MM_I$
we must have,
$$s_*(p, t)> \min (\ell_*(p,t), \de_*).$$
Furthermore, for $0\le t(p)-t\le  \min (\ell_*(p,t), \de_*)$ the foliation 
$S_t=\Si_t\cap \NN^-(p)$ is smooth. For these values of $t$ the metrics $\si_t$ on ${\Bbb S}^2$, 
obtained by restricting 
the metric $g_t$ on $\Si_t$ to $S_t$ and then pulling it back to ${\Bbb S}^2$ by the exponential
map $\GG(t,\cdot)$, and the the null lapse $\varphi^{-1}=\g(\L,\T)$ satisfy
$$
|\varphi-1|\le \eps,\qquad |\si_t(X,X)-\si_0(X,X)|\le \eps \si_t(X,X),\quad \forall X\in T{\Bbb S}^2
$$
where $\si_0$ is the standard metric on ${\Bbb S}^2$ and $\eps>0$ is a sufficiently small constant
dependent on $\de_*$.

Finally, there exists a universal constant $c>0$ such that 
$$
i_*(p)\ge c\,  \min (\ell_*(p,t), \de_*).
$$
\end{theorem}
We assume for the moment the conclusions of Theorem \ref{thm:main-conditional-t}
and proceed with the proof of  Main Theorem I.

\section{Proof of theorem I}
According to Theorem \ref{thm:main-conditional-t} the desired conclusion of  Main Theorem I
will follow after finding a small constant $\de_*$ dependent only on $\de_0, r_0,N_0, \RR_0$
and
$\KK_0$  with the property that $\ell_*(p,t)\ge \de_*$. We fix $\de_*$, to be chosen later,
and assume that
$\ell_*(p,t)<\de_*$. Recall that $g_t$ and $\si_t$ denote the restrictions of the space-time metric 
$\g$ to respectively $\Si_t$ and $S_t$, while $\si_0$ is the push-forward of the standard metric 
on ${\Bbb S}^2$ by the exponential map $\GG(t,\cdot)$. The latter is clearly well-defined 
for the values $t(p)-\ell_*(p,t)<t\le t(p)$.

We now record three statements consistent with the assumptions of the 
Main Theorem and conclusions of Theorem \ref{thm:main-conditional-t}.

{\bf A1.}\quad  There exists a constant $c_\NN=c_\NN(p)>\ell_*(p,t)$ 
such that $\NN^{-}(p)$  has no conjugate terminal points in the
 time slab $[t(p)-c_\NN,t(p)]$.

 {\bf A2.}\quad
 The  metric $\si_t$     
remains close to the metric $\si_{0}$, i.e.
 given any vector
$X$ in $T S_t$ we have,
\be{eq:small-si}
|\si_t(X,X)-\si_0(X,X)|<\ep \,\si_t(X,X)
\end{equation}
uniformly for all $t(p)-\ell_*(p,t)<t\le t(p)$. The  null lapse $\varphi=\g(\L,\T)$ also
verifies,
\be{eq:estim-varphi-ep}
|\varphi-1 | \le \ep
\end{equation}

{\bf A3.}\quad  There exists a neighborhood $\OO=I_p\times U_p$ of $p$ 
 and a system of   coordinates $x^\a$ with $x^0=t$, the time function introduced above,  
relative to which the metric $\g$ is close to the Minkowski metric $\m(p)=-n(p)^2
dt^2+\de_{ij} dx^i dx^j$,
\be{eq:spacetime-coords-small}
 |\g_{\a\b}-\m_{\a\b}(p)|<\ep
\end{equation} 
uniformly at all points in $\OO$. The set $U_p$ contains the geodesic ball $B_{r_0}(p)$ 
and $\sup_{t,t'\in I_p} |t-t'|\ge r_0$. We may assume that $r_0>>\de_*$. In particular, 
$B_{t,2(t(p)-t)}\subset \OO$ for any $t\in [t(p)-r_0/3,t(p)]$, where $B_{t,a}$ denotes 
the Euclidean ball of radius $a$ centered around the point on $\Si_t$ with the same 
coordinates $x=(x^1,x^2,x^3)$ as the point $p$.

{\bf A4.}\quad   The space-time region $\MM_I=\cup_{t\in\I}\Si_t$
contains a future, compact set
${\cal D}\subset\MM_I$  such that there exists a positive constant $\de_0$
 with the property that, for any point $q\in {\cal D}^c$, 
we have $\ell_*(q,t)>\de_0$.

{\bf Remark 1.}\quad As a first consequence of {\bf A1} we infer that 
there must exist a largest value of $t>t(p)-c_{\NN}$ with the property that 
two distinct null geodesics originating at $p$ intersect at time $t$. 
Indeed let $t_0\ge t(p)-c_\NN $
be the  supremum of  such values\footnote{Note that $t_0=t(p)-\ell_*(p,t)>t(p)-\de_*$.} 
of $t$  and    let $(q_k,
\la_k,\ga_k)$ be a sequence of points $q_k\in \NN^{-}(p)$ and distinct null geodesics  
$\la_k, \ga_k$from $p$
intersecting at $q_k$,  with  increasing $ t(q_k)\ge t(p)-c_\NN$.  By compactness we
may assume that $q_k \to q\in\Si_t$, $t(q_k)\to t=t(q)=t_0$ and  $\la_k \to \la$, $\ga_k\to
\ga$,  with both
$\la, \ga$ null geodesics passing through  $p$ and $q$. 
We claim that
$\ga\neq\la$  and that $q$ is a cut locus terminal point of $\NN^{-}(p)$. 
Indeed if $\ga\equiv \la$ then for a sequence of positive constants $0<\eps_0\to 0$ 
we could find an increasing sequence of indices $k$ such
that for null geodesics $\ga_k, \la_k$ we have that 
$$
\g(\dot\ga_k(0),\dot\la_k(0))=\eps_0,\qquad \ga_k(t(q_k))=\la_k(t(q_k)),\quad
t(q_k)>t(p)-c_{{\cal N}}.
$$
This leads to a contradiction, as by assumption the exponential map 
$\GG(t,\cdot)$ is a local diffeomorphism for all $t>t(p)-c_{{\cal N}}$.

 {\bf Remark 2.}\quad As consequence of {\bf A2} we infer that, for all  $t>
t-c_\NN$, the distances  on $S_t$ corresponding to the metrics $\si_t$ and $\si_0$ are
comparable,
\be{eq:comp-dist-on-S}
(1-\ep) d_{0}(q_1,q_2)\le d_{\si}(q_1,q_2)\le (1+\ep) d_{0}(q_1,q_2)\qquad 
\forall \, q_1,q_2\in S_t
\end{equation}
or, equivalently,
\be{eq:comp-dist-on-S2}
(1+\ep)^{-1} d_{\si}(q_1,q_2) \le d_{0}(q_1,q_2)\le (1-\ep)^{-1}d_{\si}(q_1,q_2)
\qquad \forall \, q_1,q_2\in S_t
\end{equation}
 {\bf Remark 3.}\quad Similarly, as a consequence of {\bf A3} the  distances on
$\Si_t\cap\OO$ corresponding to  the induced metric $g$ and the euclidean metric $e$ are
also comparable,
\be{eq:comp-dist-on-Si}
(1-\ep) d_{e}(q_1,q_2)\le d_{g}(q_1,q_2)\le (1+\ep) d_{e}(q_1,q_2),\qquad 
\forall \, q_1,q_2\in \Si_t\cap\OO, \,\,t\in I_p
\end{equation}
Observe also that, since $\si$ is the metric induced by $g$ 
on $S_t$,
\be{eq:compare-dist-g-si}
d_g(q_1,q_2)\le d_{\si}(q_1,q_2),\qquad \forall \, q_1,q_2\in S_t,\,\,\,t(p)-\ell_*(p)<t\le t(p)
\end{equation} 
 {\bf Remark 4.}\quad 
 In what follows we will assume, without loss of generality, that $t(p)=0$ and the 
 $x=(x^1,x^2,x^3)$ coordinates of $p$ are $x=0$. Without loss of generality
 we may also assume that $n(p)=1$. Indeed if $n(p)\neq 0$ we can rescale the 
time variable $t=\tau/n(p)$ such that relative to the new time we have
$\g=-\frac{n^2}{n(p)^2} d\tau ^2 +g_{ij} dx^i dx^j$. Once we find a  convenient value
for $\de_*'$ such that no distinct past  null geodesics from $p$ 
intersect   for values of $|\tau|\le \de_*'$ we  find the desired
  $\de_*=\de_*'\c n(p)\ge \de_*' N_0^{-1}$.

According to Remark 1 we can find a largest value of time $t_*<t(p)$ where two distinct 
null geodesics from $p$ intersect, say at point $q$ with $t(q)=t_*$. Our next result will imply
that at $q$ the angle between projections\footnote{defined relative to the decomposition 
$X=-X^0 \T + \underline X$, where $\underline X\in T_q\Si_{t_*}$.}
 of $\dot \ga_1(t_*)$ and $\dot\ga_2(t_*)$ 
onto $T_q\Si_{t_*}$ is precisely $\pi$.

\begin{lemma}\label{le:opposite}  Let $\,
t_*=t_*(p)<0$ be
 the  largest\footnote{We assume that such a point exists. } value  of $t$
 such that  that there exist two distinct
past directed null geodesics $\ga_1, \ga_2$ from $p$ intersecting at $q$ with $t(q)=t_*$. 
 Assume also that the exponential map $\GG=(t,\om)\to \ga_\om(t)$
is a global  diffeomorphism from $(t_*(p), t(p)]\times \SSS^2$
to its image on $\NN^{-}(p)$ and a local diffeomorphism in a neighborhood of $q$.
 Then,
at $q$, the projections of $\dot\ga_1(t_*)$ and $\dot\ga_2(t_*)$ onto 
$T_q\Si_{t_*}$ belong to the same line and point in the opposite directions.

\end{lemma} 

{\bf Remark 5.}\quad Similar statement also holds for future directed null geodesics.

\begin{proof}:\quad The distinct null geodesics $\ga_1$, $\ga_2$ can be identified with 
the null geodesics $\ga_{\omega_1}, \ga_{\omega_2}$ with 
$\omega_1\ne\omega_2\in {\Bbb S}^2$, two distinct directions in the tangent space 
$T_p{\cal M}$. 

By assumptions $\ga_{\om_1}(t_*)=\ga_{\om_2}(t_*)=q$ and 
there exist disjoint   neighborhoods $\VV_1$ of $(t_*, \om_1)$
and $\VV_2$ of $(t_*, \om_2)$ in ${\Bbb R}\times \SSS^2$ such that 
 the restrictions of $\GG$ to $\VV_1$, $\VV_2$ are both diffeomorphisms.
We can choose  both neighborhoods to be of the form 
$\VV_i=(t_*-\ep, t_*+\ep)\times W_i$  with $\om_i\in W_i$ for $ i=1,2$.
Let $\GG_t(\om)  = \GG(t,\om) $ and define $S_{t, i}=\GG_t( W_i), i=1,2 $.
They are both pieces of embedded 2-surfaces in $\Si_t$, $t\in (t_*-\ep_0, t_*+\ep_0)$
for some $\eps_0>0$ 
and, as the exponential map $\GG(t,\cdot)$ is assumed to be a global diffeomorphism
for any $t>t_*$,  they are 
disjoint for all $t>t_*$.

 For $t=t_*$ the surfaces $S_{t_*,i}$ intersect at the point $q$.
We claim that the tangent spaces of $T_q(S_{t_*, 1})$ and $T_q(S_{t_*, 2})$
must coincide in $T_q(\Si_*)$. Otherwise, since  $T_q(S_{t_*, 1})$ and $T_q(S_{t_*, 2})$ are two dimensional hyperplanes in a three dimensional space $T_q\Si_{t_*}$, they intersect transversally
and by an implicit function we conclude that the surfaces 
$S_{t_*,i}$ also intersect transversally at $q$. The latter is impossible, as $S_{t,1}, S_{t,2}$
are continuous families of 2-surfaces {\it disjoin}t for all $t>t_*$.

\end{proof}

The following  lemma is a consequence of 
{\bf A3} and the normalization made 
in Remark 4. Recall that ${\cal I}^-(p)$ denotes the causal past of
point
$p$, i.e.,  it consists  of all points
which can be reached by continuous, past time-like curves from $p$, and 
$\NN^-(p)$ is the boundary of  ${\cal I}^-(p)$.
\begin{lemma}\label{le:inclusions} Let $t\in [-r_0/3, 0]$ and 
 and let $p_t$ be the point on $\Si_t$  which 
has the same   coordinates $x=(x^1, x^2, x^3)=0$ as $p$.
Let $B_{t,r}=B(p_t, r)\subset\Si_t$ be  the euclidean ball centered at $p_t$ 
of radius $r$ and $B_r^c$ its complement in $\Si_t$. Then,
\be{eq:inclusions}
B_{t, (1-3\ep)|t|}\subset \II^{-}(p)\cap\Si_t, \quad B_{t, (1+3\ep)|t|}^c\cap 
\left (\II^{-}(p)\cup \NN^{-}(p)\right)=\emptyset .
\end{equation}
\end{lemma}
\begin{proof}.\quad According to Remark 4 $p$ has coordinates $t=0, x=0$ 
 and $n(p)=1$.  Hence, according to \eqref{eq:spacetime-coords-small},  $|n-1|< \ep$   
 and  $|g_{ij}-\de_{ij}|<\ep$. The point $p_t$ 
has coordinates $(t,0)$, $t>-r_0/3$. 
Let $q\in B_{t,(1-2\eps)| t|}$  of coordinate $(t,y)$ and 
$\ell(\tau)=(\tau, y \frac {\tau}t)$
 be the straight segment 
connecting $p$ with $q$.  Thus, in view of \eqref{eq:spacetime-coords-small}, 
for sufficiently small $\ep$,
\beaa
\g(\dot\ell(\tau),\dot\ell(\tau))&=& \m(\dot\ell(\tau),\dot\ell(\tau)) +
 (\g-\m)(\dot\ell(\tau),\dot\ell(\tau))\\
&\le& (-1+\frac{|y|^2}{t^2}) +\eps( 1+ \frac{|y|^2}{t^2}) =-1+\ep +(1+\ep)\frac{|y|^2}{t^2}\\
&<&-1+\ep+(1+\ep)(1-2\ep)^2=-2\ep+O(\ep)^2<0.
\eeaa
Thus  $q$ can be reached by a time-like curve from $p$, therefore
 $q\in \JJ^-(p)$.

 On the other hand, if $\ell(\tau)=(\tau,x(\tau))$ is an arbitrary  causal curve
from $p$ then, 
\beaa
0\ge \g(\dot\ell(\tau),\dot\ell(\tau))&=& \m(\dot\ell(\tau),\dot\ell(\tau)) +
(\g-\m)(\dot\ell(\tau),\dot\ell(\tau))\\
&\ge& (-1+|\dot x|^2) -\eps(1+|\dot x|^2)= -(1+\ep)+ (1-\ep)|\dot x|^2.
\eeaa
 Therefore  $|\dot x|\le\frac{ 1+\ep}{1-\ep}<1+2\ep+O(\ep^2)$ and thus,
for sufficiently small $\ep>0$,  $|x(\tau)|\le  (1+3\ep)\tau $. Consequently
points $q$ in the complement of the ball $B_{t,(1+3\ep)t}$ cannot
be reached from $p$ by a causal curve.
\end{proof}
\begin{corollary}\label{Corr:1}    Any 
 continuous  curve $x(\tau)\subset \Si_t$  between  two points 
$q_1\in  B_{t,(1-3\eps)|t|}$ and $q_2\in B_{t,r(1+3\eps)|t|}^c$
 has to intersect 
$\NN^-(p)\cap \Si_t$.
\end{corollary}
\begin{proof}:\quad  This is an immediate consequence
of  proposition \ref{le:inclusions} and the fact that
both $\II^{-}(p)$ and $\left (\II^{-}(p)\cup\NN^{-}(p)\right)^c$ are 
 connected open, disjoint, sets.
\end{proof}
{\bf Remark 6.}\quad Observe that the argument used in the 
proof of  proposition \ref{le:inclusions} also shows the inclusion,
\be{eq:another-incl}
\NN^{-}(p)\cap \Si_t\subset B_{t, (1+3\ep)|t|}\cap B^c_{t, (1-3\ep)|t|},\qquad 
\forall t:\,\, -r_0/3\le t\le 0.
\end{equation}

We are now ready prove the following, 
\begin{proposition}\label{prop:igor} Assume ${\bf A1}, {\bf A2}, {\bf A3}$ satisfied and 
$\ell_*(p,t)<\de_*<<r_0$. Then, no  
 two null geodesics originating at $p$ and also opposite\footnote{i.e.,  two null geodesics $\ga_1,\ga_2$ with the property that 
$\ga_1(0)=\ga_2(0)=p$ and 
the projections of the tangent vectors $\dot\ga_1(0), \dot\ga_2(0)$ to 
 $T_p\Si_{t(p)}$ belong to the same line and point in the opposite directions.} 
 at
 $p$ can intersect  in the slab
$[t(p)-\ell_*(p,t), t(p))$.
\end{proposition}
{\bf Remark.}\quad 
Modulo the assumption that the intersecting geodesics have to be opposite at $p$, 
Proposition \ref{prop:igor} gives the desired contradiction and implies the Main Theorem.
Indeed, Remark 1. implies that if $\ell_*(p,t)<\de_*$ then there exist two distinct null geodesics
from $p$ necessarily intersecting at time $t(p)-\ell_*(p,t)$, which contradicts the above proposition. 
The extra assumption that the geodesics are opposite at $p$ will be settled below by showing 
existence of a point $p\in \MM_I$ with the property that there exist two null 
geodesics from $p$ intersecting precisely at time $t(p)-\ell_*(p,t)$, which are also opposite at $p$.

\begin{proof}:\quad Once again we set $t(p)=x(p)=0$ and $n(p)=1$. We now argue by
contradiction.
 Assume that there exist two
opposite null geodesics 
$\ga_1\ne \ga_2$  from $p$ such that $\ga_1(t_*,\om_1)=\ga_2(t_*,\om_2)$, where $t_*$ is the 
first time of intersection of all such geodesics, with
$t_*\ge -\ell_*(p,t)$.  We choose the time $t_0>t_*$ such that the distance 
$d_g(\ga_1(t_0,\om_1),\ga_2(t_0,\om_2))< t_0 \, \eps/2$. 
Let $q_1=\ga_1(t_0,\om_1)$ and $q_2=\ga_2(t_0,\om_2)$. Note that by our assumptions 
the exponential map $\GG(t,\cdot)$ is a global diffeomorphism for all $0<t\le t_0$. Our 
assumption also implies that $\omega_1$ and $\omega_2$ represent antipodal points on
${\Bbb S}^2$.

We consider the set, 
$$
\Omega=\{\omega\in {\Bbb S}^2:\,\,d_{{\Bbb  S}^2}(\om,\om_2)<\frac \pi 4\}
=\{\omega\in {\Bbb S}^2:\,\,d_{{\Bbb  S}^2}(\om,\om_1)\ge \frac{3\pi} 4\}.
$$
Then, in view of \eqref{eq:comp-dist-on-S} and \eqref{eq:compare-dist-g-si},  the 
set $\tilde{\Om}_{t_0}={\cal G}(t_0,\Omega)$ has the
property that\footnote{Recall that $d_0$ denotes the distance function on $S_t$ defined 
with respect to the metric $\si_0$ obtained by pushing forward the standard metric on 
${\Bbb S}^2$ by the exponential map $\GG(t,\cdot)$.}, 
$$
 d_g(q_2,q)\le d_\si(q_2,q)<(1+\ep) d_0(q_2,q)=(1+\ep)|t_0|\frac \pi 4,\quad \forall q\in
\tilde{\Om}_{t_0}.
$$
Thus, in view of $d_g(q_1,q_2)< t_0 \,\eps/2$  and  the triangle inequality,
\be{eq:small}
d_g(q_1,q)\le (1+\ep)|t_0|\frac \pi 4+\frac \eps 2 t_0\le \frac\pi 4|t_0| \left (1+O(\ep)\right),
\quad \forall q\in
\tilde{\Om}_{t_0}.
\end{equation} Thus, taking into account \eqref{eq:comp-dist-on-Si},
\be{eq:small2}
d_e(q_1,q)\le\frac{1}{1-\ep}d_g(q_1,q)\le \frac\pi
4|t_0|\left (1+O(\ep)\right),\quad
\forall q\in
\tilde{\Om}_{t_0}.
\end{equation}
On the other hand, from \eqref{eq:comp-dist-on-S}, since any 
point in the complement of $\tilde{\Om}_{t_0}$ in  $\NN^-(p)\cap
\Si_{t_0}$ lies in the image by $\GG$  of the complement of $\Omega\subset\SSS^2$,
\be{eq:compl-of-tildeO}
d_\si(q_1,q)\le  (1+\eps)d_0(q_1,q)\le \frac{3\pi} 4(1+\eps)|t_0|,\quad \forall
q\in\left(\NN^-(p)\cap
\Si_{t_0}\right)\setminus \tilde{\Om}
\end{equation}
Observe that, since $q_1\in \NN^{-}(p)\cap\Si_{t_0}$ and using 
\eqref{eq:another-incl},   we have 
$q_1\in  B_{t_0,(1+3\ep)|t_0|}\cap B^c_{t_0,(1-3\ep)|t_0|}$. Thus if
$(t_0, y)$ are the coordinates of $q_1$ we must have 
$$1-3\ep\le \frac{|y|}{|t_0|}\le 1+3\ep.$$
Thus, for sufficiently small $\ep>0$,  the point 
$( t_0, -(1-7\ep\, |t_0|)y)\in B_{t_0,(1-3\ep)|t_0|}$ 
while 
$( t_0, -(1+7\ep\, |t_0|)y)\in B^c_{t_0,(1+3\ep)|t_0|}$. Let
$y(\tau)=-\big(14\ep\tau+(1-7\ep\,|t_0|)\big)\c y$, with  $\tau\in [0,1]$ and $I$
the segment
$I(\tau)=\big(t_0, y(\tau)\big)$. Observe  that all points of $I$
 are within  euclidean  distance $O(\ep\, |t_0|)$ from the point $q^{opp}=(t_0, -y)$.  
 Clearly, the extremities of $I$
verify,  $I(0)\in B_{t_0,(1-3\ep)|t_0|}$
and
$I(1)\in B^c_{t_0,(1+3\ep)|t_0|}$. To reach a
contradiction with Corollary \ref{Corr:1}   we will show that  in $I$ does not intersect
$\NN^-(p)\cap\Si_{t_0}$.

We first show that $I\cap \tilde{\Om}_{t_0}=\emptyset$. Indeed, if $q\in\tilde{\Om}_{t_0}$,
we have, from\eqref{eq:small2},
\beaa
d_e(q_1,q)\le\frac \pi 4|t_0|\left (1+O(\ep)\right)<|t_0|
\eeaa
while, if $q\in I$, 
\beaa
d_e(q_1,q)\ge |y|(2-7\ep\,|t_0|)\ge |t_0| (1-3\ep)(2-7\ep\,|t_0|)>|t_0|
\eeaa

Now, assume by contradiction, the existence of    $q\in \left(\left(\NN^-(p)\cap
\Si_{t_0}\right)\setminus
\tilde{\Om}_{t_0}\right)\cap I$. From \eqref{eq:compl-of-tildeO} we must infer that, 
$$
d_\si(q_1,q)\le \frac{3\pi} 4(1+\eps)|t_0|.
$$
On the other hand 
let $x(\tau),\,  \tau\in[0,1], $ be the $\si$-geodesic connecting $q_1$ and $q$
 in $S_{t_0}$.
Since  $S_{t_0}=\NN^-(p)\cap \Si_{t_0}$ is contained in the set 
$B_{t_0,(1+3\eps)|t_0|}\setminus B_{t_0,(1-3\eps)|t_0|}^c$   so is the entire  curve 
$x(\tau)$ for all $0\le \tau\le 1$. Now observe that the euclidean distance of any curve
which connects $q$ and $q^{opp}$ while staying outside $B_{t_0,(1-3\eps)|t_0|}$ must
be greater than $\pi(1-3\ep)|t_0|$. Since all points in $I$ are within
euclidean distance $O(\ep\,|t_0|)$ from $q^{opp}$ we infer that 
$$
\int_0^1 |\dot x(\tau)|_e d\tau \ge \pi|t_0|\left (1-O(\ep)\right).
$$
This  implies that, 
$$
d_\si(q_1,q)=\int_0^1 |\dot x(\tau)|_g d\tau \ge (1-\eps)
\int_0^1 |\dot x(\tau)|_e d\tau \ge \pi|t_0|\left (1-O(\eps)\right).
$$
 which is a contradiction. Thus
 $I$ does not intersect $\NN^{-}(p)\cap\Si_{t_0}$.

\end{proof}
The proof of proposition \ref{prop:igor} depends 
on the fact that the intersecting  null
geodesics from $p$ are opposite to each other in the tangent space
$T_p(\Si_{t(p)})$. According to lemma \ref{le:opposite} we 
know that, at the first time $t$  when two past directed
null geodesics $\ga,\la$  from $p$ intersect at a point $q$, 
they must intersect opposite
to each other. However  the situation is not entirely symmetric,
as there may exist another pair of future directed  null geodesics $\ga',\la'$
 from $q$ which intersect at a point $p_1$   with $t(p_1)$
 strictly smaller than $t(p)$. We can then repeat the procedure
with $p$ replaced by $p_1$ and with a new pair of null geodesics 
$\ga_{1},\la_{2}$ from $p_1$ intersecting  at $q_1$  with 
$t(q_1)$ the smallest value of $t$ such that any two null geodesics from
$p_1$ intersect on $\Si_t$. Proceeding by induction we can construct a
 sequence of  points $p_k$, $q_k$ with $t(p_k)$
monotonically decreasing and $t(q_k)$ monotonically increasing,
and sequence of    pairs of distinct  null geodesics
$\ga_k,\la_k$ passing through  both $p_k$ and $q_k$. Our
construction
also insures that  at $q_k$ the  geodesics $\ga_k,\la_k$
are opposite to each other. We would like to pass to limit
and  thus obtain two  null geodesics which  intersect    each other
at  two distinct points.  This procedure is behind the proof of the following
 
\begin{proposition}
Assume that the region  $\MM_I$ verifies  {\bf A4}. 
Then, if there exist two distinct  null geodesics $\la_0$, $ \ga_0$
intersecting at two points $p_0,q_0$ such that  $0< t(p_0)-t(q_0)< \de_*$,
then there must exist a pair of null geodesics $\la,\ga$ intersecting at points 
$p, q$ with $t(q_0)\le t(q)<t(p)\le t(p_0)$ which are opposite at both $p$ and $q$
\end{proposition}
\begin{proof}:\quad   Let 
\begin{equation}\label{eq:De}
\Delta t:=\min_{p,q\in \MM_I} t(p)-t(q)
\end{equation}
such that there exists a pair of distinct past directed null geodesics
originating at $p$ and intersecting at $q$. By the assumption of the proposition 
$\De t<\de_*$. On the other hand, for all 
points $p\in {\cal D}^c$, where the set ${\cal D}$ is that of the condition {\bf A4}, 
we have that $\ell_*(p,t)>\de_0$. Assuming, without loss of generality that $\de_*<\de_0$,
we see that it suffices to impose the restriction 
$p\in {\cal D}$ in \eqref{eq:De}. Since ${\cal D}$ is 
compact and the manifold $\MM_I$ is smooth we can conclude that $\De t>0$.

Let $p_n\in {\cal D}$ be a sequence of points such that 
$\ell_*(p_n,t)\to \Delta t$. Since for each $p_n$ with sufficiently large $n$ 
we have $\ell_*(p_n,t)<\de_*$ we may assume, with the help of Theorem  
\ref{thm:main-conditional-t}, that {\bf A1}--{\bf A4} are satisfied
for $p_n$.

Choosing a subsequence, if necessary, we can 
assume that $p_n\to p$. We claim that $\ell_*(p,t)=\Delta t$, i.e., 
there exists a pair of distinct past null geodesics
from $p$ intersecting at time $t(p)-\Delta t$, and that these geodesics are 
opposite to each other at $p$. First, to show existence of such geodesics we 
assume, by contradiction, that there exists an $\eps_0>0$ such that 
no two distinct geodesics from $p$ intersect at $t\ge t(p)-\De t-\eps_0$. 
Since by assumption $\De t<\de_*$ we may assume that $\NN^-(p)$ does not
contain points conjugate to $p$ in the slab $(t(p)-\De t-\eps_0,t(p))$. This implies 
that the exponential map $\GG_p(t,\cdot)$ is a global diffeomorphism for all 
$t\in (t(p)-\De t-\eps_0,t(p))$. Smooth dependence of the exponential map $\GG_q$ 
on the base point $q$ implies that there exists a small neighborhood ${\cal U}$ of $p$
such that for any $q\in{\cal U}$ the exponential map 
$\GG_q(t,\cdot)$ is a global diffeomorphism for any $t\in (t(p)-\De t-\eps_0/2,t(p))$. 
This however contradicts the existence of our sequence $p_n\to p$ since by construction
$\ell_*(p_n,t)\to \De t$. 

Therefore we may assume that there exists a pair of null geodesics
$\ga_1, \ga_2$, originating at $p$ and intersecting at a point $q$ with 
$t(q)=t(p)-\De t=t(p)-\ell_*(p,t)$. By Lemma \ref{le:opposite} the geodesics 
$\ga_1$ and $\ga_2$ are opposite at $q$. We need to show that they are also opposite at 
$p$. Consider the boundary of the causal future of $q$ -- $\NN^+(q)$. It contains a pair 
of null geodesics, the same $\ga_1$ and $\ga_2$, intersecting at $p$. Thus, either 
$t(p)$ is the first time of intersection among all distinct future directed 
null geodesics from $q$, in which 
case Remarks after Theorem \ref{thm:main-conditional-t}
and Lemma \ref{le:opposite} imply that $\ga_1$ and 
$\ga_2$ are opposite at $p$, or there exists a pair of null geodesics from $q$ intersecting
at a point $p'$ such that $t(p')<t(p)$. But then $t(p')-t(q)<\De t$ contradicting 
 the definition
of $\De t$.
\end{proof}
\section{Proof of Main  Theorem II}
We start with the following proposition.
\begin{proposition}
\label{prop:compare-metrics}
 Assume \eqref{eq:assumption-n}, \eqref{eq:assumption-T} verified.
 Then, if the initial metric $g$ on $\Si_0$ verifies  \eqref{eq:assumption-I1},
there exists a large constant $C=C( N_0, \KK_0 )$ such that, relative
to the induced transported coordinates in $I\times U$ we have,
\be{eq:assum-2-at-t}
 C^{-1}|\xi|^2\le g_{ij}(t,x) \xi^i\xi^j \le  C|\xi|^2, \qquad \forall
x\in U
\end{equation}
\end{proposition}
\begin{proof}:\quad  We fix a coordinate chart $ U$  and consider
 the transported coordinates $t,x^1,x^2,x^3$ on $I\times U$. Thus 
$
\pr_t g_{ij} = -\frac 12 n\, k_{ij}.
$
Let $X=X$ be a time-independent vector on $\M$ tangent
to $\Si_t$. Then,
$$
\pr_t g(X,X) = -\frac 12 n\, k(X,X).
$$
Clearly,
\beaa
|n k(X,X)|\le|nk|_g|X|_g^2\le \|nk(t)\|_{L^\infty} |X|_g^2
\eeaa
with $|k|_g^2=g^{ac} g^{bd}k_{ab} k_{cd}$ and $|X|_g^2=X^i X^j g_{ij}=g(X,X)$.
Therefore, since $\pr_t |X|_g^2=\pr_t g(X,X)$,
$$
-\frac 12  \|n k(t)\|_{L^\infty} |X|_g^2 \le \pr_t |X|_g^2 \le \frac 12   \|n
k(t)\|_{L^\infty}  |X|_g^2. 
$$
Thus, 
\beaa
 |X|_{g_0} e^{-\int_{t_0}^t\|nk(\tau)\|d\tau} \le |X|^2_{g_t} \le  |X|_{g_0}
e^{\int_{t_0}^t\|nk(\tau)\|d\tau}
\eeaa
from which \eqref{eq:assum-2-at-t} immediately follows. 
\end{proof}
\begin{corollary}\label{cor:comp-dist}  Let $p\in \Si_t$ in a coordinate
chart $U_t=\Si_t\cap(I\times U)$ with transported coordinates $(t,x^1,x^2,x^3)$.
Denote by $e$ the euclidean metric    on $U_t$   relative to the coordinates
$x=(x^1,x^2,x^3)$. Let $B_r^{(e)}(p)\subset U_t$ be an euclidean ball of radius $r$ 
centered at $p$. Then, for all $\rho\ge C r $, with $C=C(N_0,\KK_0)$  the constant  of
proposition \ref{prop:compare-metrics}, the  euclidean  ball   $B_r^{(e)}(p)$
 is included in the  geodesic balls $B_\rho(p)$,
relative to the metric $g_t$,
\beaa
  B_r^{(e)}(p)\subset B_{\rho}(p),\qquad \rho\ge C r.
\eeaa
\end{corollary}
\begin{proof}:\quad 
Let $q\in B_r^{(e)}(p)$ and  $\ga:[0,1]\to  B_r^{(e)}(p) $ be the line segment
between $p$ and $q$. Clearly, in view of \eqref{eq:assum-2-at-t},
\beaa
{\text dist}_{e}(p,q)=\int_0^1 \big (  e(\dot
\ga,\dot
\ga)\big )^{\frac 12}
\,d\tau\ge  C^{-1}\int_0^1 \big ( g_t( \ga(\tau)) (\dot
\ga,\dot
\ga)\big )^{\frac 12}d\tau \ge C^{-1}{\text dist}_{g_t}(p,q).
\eeaa
Thus for any $q\in B_r^{(e)}(p)$ we have
 ${\text dist}_{g_t}(p,q)\le C {\text dist}_{e}(p,q)\le Cr$.
Therefore $q$  belongs to the geodesic ball $B_{\rho}(p)$ for any $\rho\ge Cr$,
as desired.
\end{proof}
The Corollary allows us to get a lower bound for the volume radius.
We recall below the definition of volume radius on
a  general Riemannian manifold $M$.
The Corollary allows us to get a lower bound for the volume radius.
We recall below the definition of volume radius on
a  general Riemannian manifold $M$.
\begin{definition}
The volume radius $r_v(p,\rho)$ at point $p\in M$ and  scales $\le \rho$ 
is defined by,
$$
r_{vol}(p,\rho)=\inf_{r\le \rho}\frac{|B_r(p)|}{r^3}
$$
with $|B_r|$ the volume of $B_r$ relative to the metric $g$.
The volume radius $r_{vol}(M,\rho)$ of $M$ on scales $\le \rho$  is the infimum
of $r_{vol}(p,\rho)$ over all points  $p\in M$.
\end{definition}
Let $\rho_0$ be the positive number of the initial  assumption {\bf I1}.
 Thus every point $p\in \Si_t$  belongs to an euclidean ball
$B^{(e)}_{\rho_0}(p)$, relative to  local transported coordinates.
 Let $B_r(p)$ be a geodesic ball around $p$. According to Corollary  \ref{cor:comp-dist}
  for  any $a\le\min\{ \rho, r/C \}$  we must have  $B_a^{(e)}(p)\subset B_r(p)$.
Therefore, according to Proposition \ref{prop:compare-metrics},
 \beaa
 |B_r(p)|_{g_t} \ge |B_a^{(e)}(p)|_{g_t} =\int_{B^{(e)}_a(p)}  \sqrt{|g_t|}\,\, dx\ge
C^{-3/2}|B_a^{(e)}(p)|_{e}\ge C^{-3/2}a^3
\eeaa
This means that, for all $r\le C\rho$, 
\beaa
|B_r(p)|\ge C^{-3/2} (r/C)^3
\eeaa
Thus, on scales $\rho'\le C\rho$, $\rho\le \rho_0$ we must have,  $r_{vol}(p,\rho')\ge
C^{-9/2}$. Choosing $\rho\le \rho_0$ such that $C\rho=1$ we deduce the following,

\begin{proposition}
\label{prop:vol-radius}
 Under the   assumptions {\bf I1} as well as \eqref{eq:assumption-n},
 \eqref{eq:assumption-T}
 there exists a  sufficiently small
constant
  $v=v(I_0, \rho_0, N_0, \KK_0)>0$,
depending only on $I_0,\rho_0$, $N_0$,  $\KK_0$,
such that the volume radius of each $\Si_t$, for scales
 $\le 1$,  is bounded  from below,
\beaa
r_{vol}(\Si_t, 1)\ge v.
\eeaa
\end{proposition}
We rely on proposition \ref{prop:vol-radius} to prove 
the existence of good local space-time coordinates on $\MM_I$.
The key to our construction is the following general result, based
on Cheeger -Gromov convergence of Riemannian manifolds,
see  \cite  {A2} or  Theorem 5.4. in \cite{Pe}.
\begin{theorem}
\label{thm:peterson}
Given\footnote{An appropriate version
of the theorem holds in every dimension  $N$ with
an $L^p$ bound of the Riemann curvature tensor and $p>N/2$.}  $\La>0$, $v>0$ and $\ep>0$
there exists an $r_0>0$ such that on  any  3- dimensional, complete, Riemannian maniflod
$(M,g)$ with $\|R\|_{L^2}\le \La$ and  volume 
radius, at scales $\le 1$ bounded from below by $v$, i.e.,  $r_{vol}(M,1)\ge v$,
verifies the following property:

Every  geodesic ball $B_r(p)$, with $p\in M$ and $r\le r_0$
admits  a system of harmonic coordinates $x=(x^1, x^2, x^3)$ 
relative to which we have,
\bea
(1+\ep)^{-1}\de_{ij}\le g_{ij}&\le &(1+\ep) \de_{ij}\label{eq:harm-coord1}\\
r\int_{B_r(p)} |\pr^2 g_{ij}|^2 dv_g &\le&\ep  \label{eq:harm-coord2}
\eea
\end{theorem}
We  apply this theorem for the family of complete  Riemannian
 manifolds $(\Si_t, g_t)_{t\in I}$, for $p=2$. According to  proposition
\ref{prop:vol-radius} we have a uniform lower bound for the volume radius
$r_{vol}(\Si_t, 1)$. On the other hand we also have a uniform
bound on the $L^2$ norm
of the Ricci curvature tensor\footnote{which coincides with the
full Riemann curvature tensor in three dimensions. }. Indeed, according to proposition
\ref{prop:L2-curv} of the next section, there exists a constant $C=C( N_0, \KK_0)$ 
such that,  for any $t\in I$,
\beaa
\|\R(t)\|_{L^2}\le C(N_0, \RR_0) \|\R(t_0)\|_{L^2}= C \RR_0 .
\eeaa
Therefore,   for any $\ep>0$, there exists $r_0$ depending only  on
 $\ep, I_0,\rho_0, N_0, \KK_0, R_0$
such that  on any geodesic ball,  $B_{r}\subset \Si_{t}$, $r\le r_0$,
centered at a point $p_t\in \Si_{t}$,  there exist local coordinates relative to which
 the metric $g_t$ verify conditions \eqref{eq:harm-coord1}-\eqref{eq:harm-coord2}.
Starting with any such  coordinate system  $x=(x^1,x^2, x^3)$  we consider a 
 cylinder $J\times B_r$, with $J=\big(t-\de, t+\de\big)\cap I$  and the  associated
transported coordinates $(t,x)$  for which \eqref{eq:g-transported} holds, i.e.
\beaa
\g=-n^2 dt^2 + g_{ij} dx^i dx^j,
\eeaa
Integrating equation  \eqref{eq:t-sec-fund} and using assumptions \eqref{eq:assumption-n},
\eqref{eq:assumption-T} we derive, for all $t'\in J$ and $\de$ 
sufficiently small,
\beaa
|g_{ij}(t', x)-g_{ij}(t, x)|&\le 2& \int_J \|nk(s)\|_{L^\infty}ds
\le  2 N_0 |J| \sup_{t\in J}\|k(t)\|_{L^\infty} \\
&\le &2 N_0\frac{|J|}{|I|}\KK_0\le \ep
\eeaa
 provided that $4\de |I|^{-1} N_0\KK_0<\ep$.
On the other hand, according to \eqref{eq:harm-coord1}  we have 
for all $x\in B_r$,
\beaa
|g_{ij}(t,x)-\de_{ij}|\le \ep
\eeaa
Therefore, for sufficiently small  interval $J$,
whose size $2\de$ depends only on $N_0$, $\KK_0$
and $\ep>0$,
we have, for all $(t', x)\in J\times B_r$,
\bea
|g_{ij}(t',x)-\de_{ij}|\le 2\ep
\eea

On the other hand assumption \eqref{eq:assumption-T}
 also provides us with a bound for $\pr_t \log n$, i.e.    
$|I|\c\sup_{t\in I}\|\pr_t \log
n(t)\|_{L^\infty}\le \KK_0$. Hence also,
\beaa
|J|\sup_{t\in J}\|\pr_t n(t)\|_{L^\infty}\le  N_0^{-1} \frac{|I|}{|J|}\KK_0
\eeaa   Therefore, with a similar   choice
of
$|J|=2\de$ we have,
\beaa
|n(t',x)-n(t,x)|\le 2\de   N_0^{-1} \frac{|I|}{|J|}\KK_0<\ep.
\eeaa
Now, let  $n(p) $ be the value of the lapse $n$   the center  $p$ of
$B_r\subset\Si_t$.  Clearly, for
all
$x\in B_r$,
\beaa
|n(t,x)-n(p)|\le r\|\nab n\|_{L^\infty(B_r)}\le 
rN_0^{-1}\|\nab\log n\|_{L^\infty(B_r)}\le rN_0^{-1}|I|^{-1}\c \KK_0\le \ep
\eeaa
provided that  $rN_0^{-1}|I|^{-1}\c \KK_0<\ep$. Thus,
for all $(t', x)\in J\times B_r$,
\bea
|n(t',x)-n(p)|\le 2\ep
\eea
This concludes the proof of the following.
\begin{proposition}
\label{prop:cond-C}
Under assumptions {\bf I1}, {\bf I2} as well as \eqref{eq:assumption-n} and
\eqref{eq:assumption-T} the globally hyperbolic region of space-time $\MM_I$
verifies assumption {\bf C}. More precisely, for every $\ep>0$ there exists
a constant $r_0$, depending only on the fundamental constants 
$\rho_0, I_0, N_0, \KK_0, \RR_0$, such that every point $p\in \MM_I$ admits a coordinate
neighborhood $I_p\times U_p $, with each $U_p$ containing a geodesic ball $B_{r_0}(p)$
 of radius $r_0$,
 and a system of transported coordinates $(t,x)$  such that,
\eqref{eq:const-r0}, 
\eqref{eq:n-1}  and 
\eqref{eq:g-de} hold true.
\end{proposition}

The proof of theorem \ref{thm:main2}
is now an immediate consequence of Theorem \ref{thm:main1}
and proposition \ref{prop:cond-C}.

\section{Radius of conjugacy}\label{Conjugacy}
The remaining part of the paper will be devoted to the proof of Theorems 
\ref{thm:main-conditional} and \ref{thm:main-conditional-t}. As was mentioned 
before the key results on the radius of conjugacy were
 obtained\footnote{An extension of these results to null hypersurfaces 
with a vertex is part of Q.Wang's thesis, Princeton University, 2006.}  in 
\cite{KR1}-\cite{KR3}
 and here we will show how to deduce Theorems 
\ref{thm:main-conditional} and \ref{thm:main-conditional-t} from these results.

A lower bound on the radius of conjugacy in \cite{KR1}-\cite{KR3} is given by the 
following theorem. Let ${\cal L}^-(p)$ denote the union of all past directed null
geodesics from $p$. Clearly $\NN^-(p)\subset {\cal L}^-(p)$. We can extend 
the null  geodesic (potentially non-smooth) vectorfield $\L$ to ${\cal L}^-(p)$
and define $S_{s_0}={\cal L}^-(p)\cap \{s=s_0\}$ a two dimensional foliation of ${\cal L}^-(p)$
by the level surfaces of the affine parameter $s$ ($\L (s)=1$). The conjugacy radii
 of $\NN^-(p)$ and ${\cal L}^-(p)$ coincide and 
 $$ 
 \NN^-(p)\cap \left(\cup_{s\le i_*(p)} S_s\right)=\cup_{s\le i_*(p)} S_s
 $$

\begin{theorem}\label{thm:conj}
Let $\varpi>0$ be a sufficiently small universal constant and let ${\cal R}(p,s)$ 
denote the reduced curvature flux, associated with $\cup_{s'\le s} S_s$, to be 
defined below. Then there exists a large constant $C_\varpi$ such that if the
radius of conjugacy $s_*(p)\le \varpi$ then ${\cal R}(p,s_*(p))\ge C_\varpi$.
\end{theorem}

To deduce Theorems \ref{thm:main-conditional} and \ref{thm:main-conditional-t}
from Theorem \ref{thm:conj} it suffices to show that the reduced curvature flux 
${\cal R}(p,s)\le C$ for all values of $s\le \min (\ell_*(p),\de_*)$, where $\de_*$ is
allowed to depend on $N_0$, $\RR_0, \KK_0$. As we shall see below the reduced curvature
flux itself is only well defined for the values of $s<i_*(p)$. For $s_*(p)\le \min (\ell_*(p),\de_*)$ we will then show that for all $s<s_*$ we have 
the bound ${\cal R}(p,s)\le C(N_0, \RR_0,\KK_0)$ and thus by Theorem \ref{thm:conj}, in fact,
$s_*(p)> \min (\ell_*(p),\de_*)$. In the latter case, we will also show that 
${\cal R}(p,s)\le C(N_0, \RR_0,\KK_0)$ for all $s<\min (\ell_*(p),\de_*)$.

\subsection{Basic definitions and inequalities}
We start with a quick review of the Bel-Robinson 
tensor and the corresponding energy inequalities
 induced by $\T$.
The fully symmetric, traceless  and divergence free  Bel-Robinson tensor
is given by 
\be{eq:Bel-Robinson}
\Q[\R]_{\a\b\ga\de}=\R_{\a\la\ga\mu}\R_{\b\,\,\de}^{\,\la\,\,\mu}+
\dual \R_{\a\la\ga\mu}\dual\R_{\b\,\,\de}^{\,\la\,\,\mu}\,\,
\end{equation}
The curvature tensor $\R$ can be decomposed
into its  electric and magnetic parts $E, H$ as follows,
\be{eq:el-magn}
E(X,Y)=<\g(\R(X,\T)\T, Y),\qquad H(X,Y)=\g(\dual\R(X,\T)\T, Y)
\end{equation}
with $\dual\R$ the Hodge dual of $\R$. One can easily check
that $E$ and $H$ are tangent, traceless 2-tensors,   to $\Si_t$ and 
that $|\R|^2=|E|^2+|H|^2$. We easily check the formulas relative
to an orthonormal frame $e_0=T, e_1, e_2, e_3$,
\bea
\R_{abc0}&=&-\in_{abs}H_{sc},\qquad \dual
\R_{abc0}=\in_{abs}E_{sc}\label{eq:R-coeff}\\
\R_{abcd}&=&\in_{abs}\in_{cdt}E_{st},\qquad 
\dual\R_{abcd}=-\in_{abs}\in_{cdt}H_{st}\nn
\eea
Observe that,
\be{eq:estim-Q}
|\Q|\le 4( |E|^2+|H|^2)
\end{equation}
and, 
\be{eq:formulas-Q}
\Q_{0000}=|E|^2+|H|^2
\end{equation}
Let $\P_\a=\Q[\R]_{\a\b\ga\de} \T^\b \T^\ga \T^\de$.
By a straightforward calculation,
\bea\label{eq:Div-P}
\D^\a \P_\a=
\frac{3}{2}\,\piT^{\a\b} \Q_{\a\b\ga\de} \T^\ga \T^\de
\eea
Therefore,
integrating in  a slab  $\MM_J=\cup_{t\in J}\Si_t$,    $J=[t_0, t]\subset I$, 
we derive the following. 
\bea
\int_{\Si_t}
 \Q_{0000}=\int_{\Si_0}  \Q_{0000}+
\frac{3}{2}\int_{t_0}^t\int_{\Si_{t'}} n \,\piT^{\a\b} \Q_{\a\b
00}\, dv_g\label{eq:first-en}
\eea
with $dv_g$ denoting the volume element on $\Si_t$.
Now,
\beaa
\big|\int_{t_0}^t\int_{\Si_{t'}} n \,\piT^{\a\b}\Q_{\a\b
00} dv_g\big|&\les& N_0\int_{t_0}^t\int_{\Si_{t'}} | \piT|(|E|^2+|H|^2)  dv_g\\
&\les& N_0\int_{t_0}^t\|\,\piT(t')\|_{L^\infty}\big(\|E(t')\|_{L^2}+\|H(t')\|_{L^2}\big) dt'
\eeaa
Thus, if we denote 
$$\QQ(t)=\int_{\Si_t} \Q_{0000}=\int_{\Si_t}(|E|^2+|H|^2)dv_g, $$
we deduce,
\beaa
\QQ(t)-\QQ(t_0)\les N_0\int_{t_0}^t
\|\pi(t')\|_{L^\infty}\QQ(t')dt'
\eeaa
and by Gronwall,
\beaa
\QQ(t)\les\QQ(t_0)\expp \big(\int_{t_0}^tN_0 \|\,\piT(t')\|_{L^\infty}dt'\big)
\eeaa
Thus, in view of \eqref{eq:assumption-T},  
\beaa
\QQ(t)\les\QQ(t_0)\, \expp{\big( N_0 \KK_0\big)}
\eeaa
We have just proved the following,

\begin{proposition}\label{prop:L2-curv} Assume that the assumptions
\eqref{eq:assumption-n} and \eqref{eq:assumption-T}
 are  true.  There exists a constant $C=C( N_0, \KK_0)$ 
such that,  for any $t\in I$,
\be{eq:L2-curv}
\|\R(t)\|_{L^2}\le C \|\R(t_0)\|_{L^2}= C \RR_0 .
\end{equation}
\end{proposition}
Instead of integrating \eqref{eq:Div-P} in the slab $\MM_J$
we  will now integrate  it in the region $\DD^{-}_J(P)=\JJ^{-}(p)\cap \MM_J$
 whose boundary consists of the null part  $\NN^{-}(p)$
 and spacelike part $  D_0(p)=\JJ^{-}(p)\cap \Si_0.$ 
We recall that $\NN^{-}(p)$ is a Lipschitz manifold
and the set of its terminal points $\TT^{-}(p)$ has measure zero relative to 
$dA_{\NN^{-}(p)}$.

Let $(\P^*)_{a\b\ga}=\in_{\a\b\ga\mu} \P^\mu$
and the associated differential form,
$^*\P=(^*\P)_{\a\b\ga} dx^\b dx^\ga dx^\de$.
We can rewrite equation \eqref{eq:Div-P}
in the form,
$d^*\P=-^* \F,$
with $(^* \F)_{\a\b\ga\de}=\in_{\a\b\ga\de} \F$,
and,
$$\F=\frac{3}{2}\, \piT^{\a\b} \Q_{\a\b\ga\de}\T^\ga \T^\de.$$
Integrating  the last expression 
 in the space-time  region $\DD^{-}_I(p)=\JJ^{-}(p)\cap \MM_J$,
with $J=[t_0, t]$,  $p\in \MM_J$,
and applying Stokes theorem
we derive,
\bea
\int_{\DD_J^{-}(p)} \dual \F&=&-\int_{\NN^{-}(p)\cap\MM_J}\dual
\P=\FF_p(\NN^{-}(p)\cap\MM_J)-\En(D_0(p))\label{eq:Flux-En}
\eea
where
\bea
\En(D_0(p))&=& -\int_{D_0(p)}\dual\P=\int_{D_0(p)}\Q(\T,\T,\T,\T) dv_g\label{eq:energy-restr}\\
\FF_p(\NN^{-}(p)\cap\MM_J)&=&-\int_{\NN^{-}(p)\cap\MM_J}\dual\P \label{eq:flux-restr}
\eea
 The energy integral \eqref{eq:energy-restr}  through $D_0(p)\subset\Si_0$
 can clearly be bounded by $\|\R(t_0)\|_{L^2}$. Moreover,
in view of proposition \ref{prop:L2-curv}
the integral $\int_{\DD_J^{-}(p)} \dual \F$ can be bounded by 
$C(\KK_0)\c \RR_0$. Therefore,
\be{eq:bound-flux}
\FF_p(\NN^{-}(p)\cap\MM_J)\les C(\KK_0)\c \RR_0
\end{equation}

We recall that the null boundary $\NN^{-}(p)$ is 
a Lipschitz manifold.
This means that  every 
point $p\in\NN^{-1}(p)$ has   a local  
 coordinate chart $U_p$  together with
local coordinate  $x^\a=x^\a(\tau,\om^1,\om^2)$
which are Lipschitz continuous. The coordinates are
such that for all fixed\footnote{ The statements here are understood  to be 
true with
the possible exception of a set of measure zero relative
to the measure  $dA_{\NN^{-}(p)}$ along $\NN^{-}(p)$ introduced just
before definition \ref{def:rad-inj}.} $\om=(\om^1, \om^2)$ the
curves
$\tau\to x^\a(\tau,\om)$ are null and for  any fixed  value $\tau$  the $2$ dimensional
surfaces 
$S_\tau$,
 given by  $x^\a=x^a(\tau,\om)$,  are space-like.
 In particular there is  a 
well defined null normal  $\frac{dx^\a}{d\tau}=l^\a$ at all points  of $U_p$
with the possible exception of a set of measure zero. Moreover we can choose our coordinate
charts such that at each point where the normal  $l$ is defined
we have $\g(l,\T)>0$, i.e. $l$ is past oriented. Observe that on such coordinate chart
$U$ we have,
\beaa
\int_U\dual \P&=&\int_U \dual \P_{\a\b\ga} dx^\a dx^\b dx^\ga
=\int_U \g(\P, l)d\tau dA_{\tau}=\int_U \Q(\T,\T,\T,l)d\tau dA_{\tau}
\eeaa
with $dA_\tau$ the volume element of the space-like  surfaces $S_\tau$.
Since $\T$ is future time-like and $l$ is null past directed  we have
$ \Q(\T,\T,\T,l)<0$. Consequently, for every coordinate chart 
$U\subset\NN^{-}(p)$, $\FF^{-}_p(U)\ge 0$, where
\be{eq:def-curv-flux}
\FF_p^{-}(U)=-\int_U \dual \P
\end{equation}
Using a partition of unity it follows that
$\FF_p(U)\ge 0$ for any $U\subset \NN^{-}(p)$ and therefore
$\FF^{-}_p(U_1)\le \FF^{-}_p(U_2)$ whenever $U_1\subset U_2\subset \NN^{-}(p)$.
We can thus identify $\FF_p(U)$ as the flux of curvature through $U\subset\NN^{-}(p)$.

Therefore we have the following:
\begin{proposition}
\label{prop:flux}Under assumptions \eqref{eq:assumption-n},\eqref{eq:assumption-T} and 
\eqref{eq:initial-L2-curv}  the  flux of curvature  in $\MM_I\cap \NN^{-}(p) $, 
  $\FF^{-}_p( \MM_I)=\FF^{-}_p( \MM_I\cap\NN^{-}(p))$,
can be bounded by a uniform constant 
independent of $p$. More precisely,
\beaa
 \FF_p( \MM_I)\le C(N_0, \KK_0) \c \RR_0.
\eeaa
\end{proposition}
\subsection{Reduced curvature flux}
\label{sect:geod-curv-flux}
Let $S_s$
be the $2$ dimensional space-like  surface of a constant affine parameter
$s$, defined by the condition $\L(s)=1$ and $s(p)=0$.  
Clearly for $s\le \de<i_*(p)$ the union of $S_s$ defines a regular foliation 
of 
$$
\NN^{-}(p,\, \de)=\cup_{s<\de} S_s
$$
At any point of  $\NN^{-}(p,\, \de)\setminus \{ p\}$ we can define a 
conjugate null vector $\LLb$ with $\g(\L, \LLb)=-2$
and such that  $\LLb$ is orthogonal to the leafs   $S_s$. 
In addition we can choose  $(e_a)_{a=1,2}$ tangent  $S_s$
such that together with $\L$ and $\LLb$ we obtain a null frame,
\begin{align}
&g(\L,\LLb)=-2,\qquad \g(\L,\L)=\g(\LLb,\LLb)=0,\nn\\
&\g(\L,e_a)=\g(\LLb,e_a)=0,\qquad \g(e_a,e_b)=\de_{ab}.\label{eq:null frame}
\end{align}
We denote by $\si$ the  restriction 
of $\g$ to $S_{s}$.
Endowed with this  metric  $S_s$ is a $2$ dimensional 
 compact riemannian manifold
 with  $\ga(e_a,e_b)=\de_{ab}$.
 Let $|S_s|$ denotes the area of $S_s$
and define 
$r=r(s)$
 by the formula
\be{eq:define-r}
4\pi r^2=|S_s|.
\end{equation}
 We say that 
a tensor $\pi$ along $\NN^{-}(p)$ is $S_s$ tangent, or simply $S$-tangent,
if, at every point of  $\NN^-(p)$   it is orthogonal to both null
vectors  $\L$ and $\LLb$. Given such a tensor, say  $\pi_{ab}$, we denote by  $|\pi|$
its length relative to the metric $\ga$, i.e. $|\pi|^2=\sum_{a,b=1}^2|\pi_{ab}|^2.$

We  denote by $\nab$  the restriction of $\D$ to $S_s$,
Clearly, for all $X, Y\in T(S_s)$,
\be{eq:def-nab}\nab_X Y=\D_X Y+\f12 <\D_X Y, \LLb>\L+\f12 <\D_X Y, \L>\LLb
\end{equation}
Given an $S-$ tangent tensor $\pi$ we define $(\nab_L \pi)$ to be the projection
to $S_s$ of $\D_\L\pi$. We write $\nabb \pi=(\nab \pi, \nab_L\pi)$
and 
$$|\nabb\pi|^2=|\nab_L\pi|^2+|\nab\pi|^2.$$
We  recall the definition of the  null second fundamental form $\chi, \chib,$
 and torsion $\ze$ associated to the
  $S_s$ foliation.
\be{eq:null lapse-nullsecondf}
 \chi_{ab}=\g(\D_{e_a}\L, e_b),\qquad \chib_{ab}=g(\D_{e_a}\LLb, e_b)\qquad 
\ze_a=\f12 \g(\D_a\L,\LLb)
\end{equation}
We also introduce,
\be{eq:null-lapse}
\varphi^{-1}=g(\T,\L),\qquad \psi_a=\g(e_a, \T)
\end{equation}
Observe that $\varphi>0$  with $\varphi(p)=1$. Also
\be{eq:dtds}
\frac{dt}{ds}=-n^{-1}\varphi^{-1}
\end{equation}
with $n$ the lapse function
of the $t$ foliation.
We now  recall the standard  null decomposition of the Riemann curvature tensor
 relative to  the $S_{s}$ foliation:
\bea
\a_{ab}&=&\R_{\L a\L b}\,,
\quad \b_a=\f12 \R_{a \L\LLb \L} ,\quad
\rho=\frac{1}{4}\R_{\LLb \L \LLb \L}\,\quad\nn\\
\quad
\si&=&\frac{1}{4}\, ^{\star} \R_{\LLb \L\LLb \L},\quad
\bb_a=\f12R_{a\LLb\LLb \L},\quad \aa_{ab} =\R_{\LLb a\LLb b}\
\eea
We  can  write  the flux along 
$\NN^-(q,\de)$,   $\de<i_*(q)$  as follows.
\beaa
\FF(p,\de)=\int_{\NN^-(q,\de)}\Q(\T,\T,\T, \L)=\int_{0}^\de ds
\int_{S_s}\Q(\T,\T,\T,
\L)dA_s
\eeaa
Observe that $ds dA_s$ is precisely the measure $dA_{\NN^{-}(p)}$.
More generally we shall use the following notation.
\begin{definition} \label{def:integr-Nqde}
Given a scalar function $f$
on $\NN^-(p,\de)$, $\de \le i_*(p)$ we denote  its integral 
on   $\NN^-(p,\de)$  to be,
 $$\int_{\NN^-(p,\de)} f=\int_{0}^\de ds \int_{S_s} f dA_s=
\int_{\NN^-(p,\de)} f \,\,  dA_{\NN^{-}(p)} .$$
Or, relative to the normal  coordinates $(s, \om)$ in the tangent space to
$p$, 
\beaa
\int_{\NN^-(p,\de)} f=\int_{0}^\de\int_{|\om|=1} f(s,\om)\sqrt{|\si(s,\om)|} ds d\om
\eeaa
where $|\si(s,\om)|$ is  the determinant of the components of the 
induced metric $\si$ on $S_s$ relative to the coordinates $s, \om$.
\end{definition}
To express the density $\Q(\T,\T,\T, \L)$
 in terms of the null components 
$\a,\b,\rho,\si, \bb$   we need to relate $\T$ to 
the null frame $\L,\LLb, e_a$. To  do this we  first introduce  
another null frame attached to the $t$ foliation.
 More precisely, at some point $q\in\NN^{-}(p,\de)$,  we let $S_t=\Si_t\cap \NN^{-}(p)$ 
 for $t=t(p)$. We define $\LLb'$ to be the null pair conjugate to $\L$
relative to $S_t$. More precisely $\g(\L, \LLb')=-2$
 and $\LLb'$ is orthogonal to $S_t$.  We complete $\L,\LLb'$
 to a full null frame on
$S_t$ by
\beaa
e_a'=e_a-\varphi \psi_a \L
\eeaa
We also have,
\beaa
\LLb'=\LLb-2\varphi\psi_a e_a+ 2\varphi^2|\psi|^2 \L
\eeaa
Now,
\beaa
\T&=&-\f12 (\varphi \L+\varphi^{-1}\LLb')=-\f12 \varphi
\L-\f12\varphi^{-1}(\LLb-2\varphi\psi_a e_a+ 2\vphi^2|\psi|^2 \L)
\eeaa
Therefore,
\be{eq:T-geod-foliation}
\T=\varphi(-\f12-|\psi|^2)\L-\f12 \varphi^{-1}\LLb+\psi_a e_a
\end{equation}
which we rewrite in the form,
\bea
\T&=&T_0+X ,\qquad T_0=-\f12\L-\f12 \LLb\\
X&=&\big(-\f12  (\varphi-1)  -\varphi |\psi|^2 \big)\L-\f12( \varphi^{-1}-1)\LLb+\psi_a e_a
\label{eq:T-geod-foliation}
\eea
Now,
\beaa
\Q(\T,\T,\T, \L)&=&\Q(T_0+X, T_0+X, T_0+X, \L)=\Q(T_0, T_0, T_0, \L)+ \Qr\\
\Qr&=& \Q(X, T_0,T_0,\L)+\Q(X, X,T_0,\L)+\Q(X, X,X,\L)
\eeaa
By a straightforward calculation,
\beaa
\Q(T_0, T_0, T_0, \L)&=& \frac 1 4 |\a|^2+
 \frac 3 2  |\b|^2+\frac 3 2  ({\rho}^2 + {\si}^2) 
+\f12 |\bb|^2
\eeaa
For $\de<i_*(p)$ we introduce the reduced flux, or geodesic
curvature flux,   
\be{eq:rondR}
\RR(p,\de)=\big(\int_{0}^\de \int_{S_s}
\big(|\a|^2+|\b|^2+|\rho|^2+|\si|^2+|\bb|^2\big)dA_s ds \big)^{1/2}
\end{equation}

On the other hand the following result can be 
easily seen from \eqref{eq:T-geod-foliation} .
\begin{lemma}
\label{le:boots-vphi-psi}
Assume that the following estimates hold on 
$\NN^-(p,\de)$, for some $\de<i_*(p)$,
\be{eq:boots-vphi-psi}
|\varphi-1|+|\psi|\le 10^{-2}
\end{equation}
Then on $\NN^-(p,\de)$,
\beaa
\Q(\T,\T,\T, \L)&\ge&  \f12 Q(T_0,T_0,T_0,
\L)\ge
\frac{1}{8}\big(|\a|^2+|\b|^2+|\rho|^2+|\si|^2+|\bb|^2\big)
\eeaa
\end{lemma}
{\bf Remark.}
We can guarantee the existence of such $\de>0$, as the initial conditions 
for $\vphi$ and $\psi$ are
$\vphi(p)=1$ and
$\psi(p)=0$. The challenge will be to extend estimate \eqref{eq:boots-vphi-psi}
to a larger region.

As an application of proposition \eqref{prop:flux}
and lemma \ref{le:boots-vphi-psi} above
we derive,
\begin{corollary}
\label{prop:reduced-flux} Let $p\in \MM_J$ and 
assume that the  estimate \eqref{eq:boots-vphi-psi}
holds on $\NN^-(p,\de)$ for some 
$\de<i_*(p)$. Then   the  reduced curvature flux  $\RR(p,\de)$
can be bounded from above  by a constant which depends only on 
$N_0$,  $\KK_0$ and the initial data  $\RR_0$. 
\end{corollary}
In view of Theorem \ref{thm:conj} and Corollary \ref{prop:reduced-flux} 
to finish the proof of Theorem \ref{thm:main-conditional} 
we need to show that there exists a constant $\de_*=\de_*(N_0, \KK_0,\RR_0)$
such that the bounds  \eqref{eq:boots-vphi-psi} can be extended to all values 
values of $s\le \min (\de_*,i_*(p))$.

We first  state  a theorem  which is an extension of Theorem \ref{thm:conj} 
and another consequence  
of  the results proved in \cite{KR1}-\cite{KR3}. We will then show simultaneously
that for all values of $s\le\min (\de_*,i_*(p))$ the reduced curvature flux 
$\RR(p,s)\le C(\RR_0,\KK_0)$ and the estimates  \eqref{eq:boots-vphi-psi}
hold true.

\begin{theorem}
\label{thm:radius of conj} Let $p\in\MM_I$ fixed and 
assume that the reduced curvature flux  verifies  $\RR(p,\de)\le C$
for some $\de\le i_*(p)$ and a positive constant $C$. 
 Let  $\varep_0>0$ be  a fixed small constant. 
 Then for all $s\le \min(\varpi,\de)$, where $\varpi$ is a small constant
 dependent only on $\eps_0$ and $C$, we have
 \be{eq:bound-chi}
|\trch-\frac{2}{s}|\le \varep_0,\qquad
\int_{0}^{s}|\chih|^2 ds'\le \varep_0.
\end{equation}
\end{theorem}
\subsection{Bounds for $\vphi$ and $\psi$}
The proof of the bounds for the reduced curvature flux and $\vphi$ and $\psi$ 
depends, in addition to 
the results stated in Theorem \ref{thm:radius of conj}
and Corollary \ref{prop:reduced-flux} on the following,
\begin{proposition}
\label{prop:estim-phi-psi} 
Let $\de_*$ be a small constant dependent only on $N_0,\KK_0, \RR_0$.  
Assume that  $\trch$ and $\chih$ 
 verify \eqref{eq:bound-chi}  for all $0\le s\le \de<\min (\de_*,i_*(p))$. Assume also
that the condition \eqref{eq:boots-vphi-psi} holds true
for $0\le s\le \de$ and let $\eps_0<10^{-1}$ in Theorem \ref{thm:radius of conj}. Then
 the following better estimate holds for all $0\le s\le \de$,
$$|\vphi-1|+|\psi|\le 10^{-3}.$$
\end{proposition} 
{\bf Remark.}\quad
The above proposition, Corollary  \ref{prop:reduced-flux} and a simple continuity argument 
allow us to get the desired conclusion that the reduced curvature flux $\RR(p,\de)$ is 
bounded by $C(N_0,\KK_0,\RR_0)$ for all $\de<\min (i_*(p),\de_*)$, which in turn, by Theorem 
\ref{thm:conj}, implies that
$s_*(p)>\min(\ell_*(p),\de_*)$.

\begin{proof}:\quad
We 
shall use the frame  $\L, \LLb', e_a'$ 
attached to the foliation $S_t$. Recall that,
\beaa
e_a'=e_a-\varphi \psi_a \L,\qquad \LLb'=\LLb-2\varphi\psi_a e_a+ 2\varphi^2|\psi|^2 \L
\eeaa
and 
$$\T=\varphi(-\f12-|\psi|^2)\L-\f12 \vphi^{-1}\LLb+\psi_a e_a.$$
We denote by $N$ the vector,
$$N=-\f12 (\varphi \L-\varphi^{-1}\LLb').$$
Observe that $\g(N,\T)=0$ while $\g(N,N)=1$. Thus $N$ is the unit
normal to $S_t$ along the hypersurface $\Si_t$.
 We can now  decompose $\L$ and $\LLb'$ as follows,
\bea
\L=-\f12 \vphi^{-1}(\T+N),\quad \LLb'=-\f12 \vphi(\T-N)
\eea

We  shall next derive  transport equations for $\vphi$ and 
 $\psi_a$ along $\NN^{-}(p)$.  
 We start with $\varphi$ and, recalling the definition
of $\piT$, we derive,
\beaa
\frac{d}{ds}\vphi&=&\frac d{ds}\,\g(\T, \L)=\g(\D_\L\T, \L)=\f12\, \piT_{\L\L}\\
&=&\frac{1}{4}\vphi^{-2}( \,\piT_{\T N}+\f12\, \piT_{NN})
\eeaa
On the other hand, writing
$$\D_\L e_a=\ddd_\L e_a-\zeta_a \L,\qquad \ze_a=\f12 \g(\D_a\L,\LLb).$$
we have,
\beaa
\nab_\L\psi_a&=&\g(\D_\L \T,e_a)- \g(\T, \L)\ze_a
\eeaa
Observe that
$
\piT_{\L e_a}=\g(\D_\L \T,e_a)+\g(\D_a\L, \T).
$
Therefore, since $\T=\varphi(-\f12 -|\varphi^2|\psi|^2)\L-\f12 \vphi^{-1}\LLb+\psi_b e_b$,
\beaa
\nab_L\psi_a&=& \piT_{\L e_a}-\g(\D_a\L, T)-\vphi^{-1}\ze_a\\
&=& \,\piT_{\L e_a}-\vphi^{-1}\ze_a-\g(\D_a\L, -\f12\vphi^{-1} \LLb+\psi_b e_b)\\
&=&\pi_{\L e_a}-\chi_{ab} \psi_b=
\,\piT_{\L e_a'}-\vphi\psi_a\pi_{\L\L} -\chi_{ab} \psi_b\\
&=&-\chi_{ab} \psi_b-\f12 \vphi^{-1}\psi_a(\,^\T \pi_{TN}+\f12\,\piT_{NN})
-\f12 \vphi^{-1}(\, \piT_{0 a'}+\, \piT_{Na'})
\eeaa
Thus the scalar $\vphi$ and the $S_s$-tangent vectorfield $\psi_a$ satisfy the equations:
\begin{align}
&\frac{d}{ds}\vphi=\frac{1}{4}\vphi^{-2}(\,^\T \pi_{TN}+\f12\,\piT_{NN}),
\label{eq:transp-vphi},\\
&\nab_\L\psi_a+\chi_{ab}\psi_b=-\f12 \vphi^{-1}\psi_a( \,\piT_{TN}+\f12\,\piT_{NN})
-\f12 \vphi^{-1}(\,\piT_{T a'}+\, \piT_{Na'})
\label{eq:transp-psi}
\end{align}
with initial conditions $\vphi(0)=1$ and $\psi(0)=0$.
In view of our main assumptions   \eqref{eq:assumption-n}, \eqref{eq:assumption-T} we  have the
obvious bounds,
\beaa
|\, \piT_{TN}|+|\,\piT_{T a'}|+|\,\piT_{NN}|+|\,\piT_{Na'}|&\les&|I|^{-1} \KK_0
\eeaa
In view of these bounds we find by integrating equation 
\eqref{eq:transp-vphi},
\beaa
|\vphi(s)-1|\les \KK_0 s/|I|
\eeaa
To estimate $\psi$ we first rewrite equation
\eqref{eq:transp-psi} in the form,
\beaa
\big |\frac d{ds} s^2 |\psi|^2 \big|\les|I|^{-1} \KK_0 s^2 |\psi|(1+|\psi)  +
s^2(|\trch-\frac{2}{s}|+
|\chih| ) |\psi|^2,
\eeaa
Integrating and using  the bounds \eqref{eq:bound-chi} for $\chih$ and $\trch$ we  obtain 
$$
|\psi(s)|^2 \les \KK_0s/|I|+\varep  s^{1/2}
$$
for any $0\le s\le \de$.  Therefore,  for $\varep\le 10^{-1}$
the desired bounds for
$\vphi$ and
$\psi$ of proposition \ref{prop:estim-phi-psi} 
 can be obtained in any interval $[0,\de]$  as
long as
$\de\c\KK_0/|I|+ 10^{-1}\de^{1/2}< <10^{-3}$. 
\end{proof}
\subsection{Proof of Theorem \ref{thm:main-conditional-t}}
The proof of the corresponding result for the $S_t$-foliation proceeds along 
the same lines as the one above for the geodesic foliation $S_s$. The connection
between the two foliations is given by the relations
\begin{align*}
&\frac{dt}{ds}=-n^{-1}\vphi^{-1},\\
&\chi_{a'b'}=\chi_{ab},\qquad \ze_{a'}=\ze_{a} - \vphi \psi_b \chi_{ab}
\end{align*}
We leave the remaining details to the reader.

\end{document}